\documentclass[leqno,11pt]{amsart}
\usepackage{amsfonts,amssymb,graphicx}
\usepackage{amsfonts}
\usepackage{amsmath}
\usepackage{amstext}
\usepackage{amsthm}
\usepackage{xspace}

\newtheorem{thm}{Theorem}[section]

\newtheorem{lem}[thm]{Lemma}

\theoremstyle{definition}
\newtheorem{defn}[thm]{Definition}
\theoremstyle{definition} \numberwithin{equation}{section}

\newcommand{\set}[1]{\left\{#1\right\}}
\newcommand{\R}{\mathbb R}
\newcommand{\Z}{\mathbb Z}
\newcommand{\Q}{\mathbb Q}
\newcommand{\N}{\mathbb N}

\newcommand{\La}{\mathbf \Lambda}
\newcommand{\eps}{\epsilon}
\newcommand{\lo}{\log_{_2}}

\newcommand{\A}{\mathcal{A}}

\newcommand{\C}{\mathcal{C}}
\newcommand{\D}{\mathcal{D}}
\newcommand{\U}{\mathcal{U}}
\newcommand{\V}{\mathcal{V}}
\newcommand{\F}{\mathcal{F}}

\newcommand{\virg}[1]{``#1"}

\newcommand{\card}{\textrm{card}}

\begin{document}
\title[Complexity for extended dynamical systems]{Complexity for
  extended dynamical systems}
\author{Claudio Bonanno}
\address{Dipartimento di Matematica Applicata, Universit\`a di Pisa,
  via F.Buonarroti 1/c, 56127 Pisa, Italy}
\email{bonanno@mail.dm.unipi.it}
\author{Pierre Collet}
\address{Centre de Physique Th\'eorique, \'Ecole Polytechnique, CNRS
  UMR 7644, F-91128 Palaiseau Cedex, France}
\email{collet@cpht.polytechnique.fr}

\thanks{{\it Acknowledgements.} The first named
author would like to acknowledge support and kind hospitality by
the Centre de Physique Th\'eorique during his stay at the \'Ecole
Polytechnique, Palaiseau Cedex, France.}%

\begin{abstract}
We consider dynamical systems for which the  spatial extension
plays an important role. For these systems, the notions of
attractor, $\eps$-entropy and topological entropy per unit time
and volume have been introduced previously. In this paper we use
the notion of Kolmogorov complexity to introduce, for extended
dynamical systems, a notion of complexity per unit time and volume
which plays the same role as the  metric entropy for classical
dynamical systems. We introduce this notion as an almost sure
limit on orbits of the system. Moreover we prove a kind of
variational principle for this complexity.
\end{abstract}
\maketitle

\section{Introduction} \label{intro}

Dynamical systems are called \virg{extended} when the spatial
extension plays an important role. They occur for example, in
nonlinear partial differential equations of parabolic or
hyperbolic type when the size of the domain is much larger than
the typical size of the structures developed by the solutions. As
in Statistical Mechanics, one can try to use the infinite volume
limit as an approximation.

For several classes of such systems, motivated by physical models, it
has been shown that one can define the semi flow of evolution in
unbounded domains acting on bounded functions with some regularity
(see for example \cite{c}, \cite{gv}, \cite{m}, \cite{C.E.}). This is
particularly convenient when studying traveling solutions or waves,
since one would not like to fix some particular boundary conditions
which restrict the nature of the solution (for example fixing a
particular spatial period). Once the dynamics has been defined in
unbounded domain, one can ask for a notion of attractor. Such a notion
was introduced by Feireisl (see \cite{Feireisl} and \cite{ms}) by
observing the system in bounded windows and inferring the result for
the unbounded domain. When the evolution equation does not depend
explicitly on space (homogeneous system), the attractor is translation
invariant and often non compact of infinite dimension. However, if
restricted to a finite window it is often a compact set. A situation
which occurs in several examples is that the functions in the
attractor are analytic and bounded in a strip around the real domain
(see for example \cite{C1.}, \cite{T.B.D.H.T.}). Compactness follows
in bounded (real) regions when using $C^{k}$ norms for example. For
such systems with non compact translation invariant infinite
dimensional attractors, one can try to define extensive quantities as
in statistical mechanics. A notion of dimension per unit volume can be
defined from the $\eps$-entropy per unit volume of Kolmogorov (see
\cite{kt}) where it was used in particular to quantify the fact that
some function spaces are larger than others. Looking for example at an
attractor composed of functions analytic in a strip and of infinite
dimension, since an analytic function is completely determined by its
data in a finite domain, the dimension observed in any finite window
will always be infinite . To avoid this uninteresting result, one
first fixes a precision $\eps>0$. One then counts for example the
minimal number $N^\eps_\Lambda$ of balls of radius $\eps$ needed to
cover the attractor in the finite window $\Lambda$. The next step is
to prove that $H(\eps)=\lim_{|\Lambda|\to\infty}\lo
N^\eps_\Lambda/|\Lambda|$ exists, and then to consider the quantity
$H(\eps)/\lo\eps^{-1}$ for small $\eps$. As mentioned above, the order
in which the limits in $\Lambda$ and $\eps$ are taken is important. If
for fixed $\Lambda$ one lets first $\eps$ tend to zero, the result is
in general infinite, while in the other order, one can get finite
results.  These ideas were applied to the attractors of various
extended systems (see \cite{C.E.1.}, \cite{E.M.Z.} and \cite{c2}).
These ideas can also be adapted to give a definition of the
topological entropy per unit volume (see \cite{C.E.2.}, \cite{C.E.3.},
\cite{c2} and \cite{zelik}). One first fixes a finite precision,
considers the maximal number $N^\eps_\Lambda(T)$ of different
trajectories one can observe in a finite window $\Lambda$ on the time
interval $[0,T]$ at this given precision. One then considers the
limits $$
h_{top}=\lim_{\eps \searrow
  0}\lim_{|\Lambda|\to\infty}\frac{1}{|\Lambda|}
\lim_{T\to\infty}\frac{\lo N^\eps_\Lambda(T)}{T}\;.  $$
Here again the
order of the limits is crucial, otherwise one may get an infinite
quantity.

Regarding similar approaches, angular limits have been proposed in
\cite{mi} for cellular automata, and for entropies along subspaces
we refer to \cite{mz} and references therein. In \cite{C.E.2.} a
similar definition was proposed for the metric entropy per unit
time and volume, however this definition involves several limits
which are up to now not known to exist.

In order to circumvent this difficulty we deal in the present
paper with the Kolmogorov complexity. For dynamical systems on a
compact phase space with an ergodic invariant measure, it is known
that the complexity per unit time of a typical trajectory is equal
to the metric entropy (see \cite{brudno}, \cite{white}). A first
advantage of the complexity per unit time is that it can be
defined for individual trajectories with initial conditions on a
full measure set. We will also see below that the complexity
satisfies some useful sub-additivity properties allowing to define
a complexity per unit time and unit volume. The strategy is the
same as for the topological entropy. We first fix a precision
$\eps$.  We then consider the complexity per unit time of a coding
of these trajectories in the window $\Lambda$ using a covering by
balls of radius at most $\eps$. We then show that this quantity
grows like the volume and define a complexity per unit time and
unit volume at a fixed precision, finally letting the precision
become infinite.

We will deal in the present paper with systems satisfying some
hypothesis inspired by the known results on extended systems (for
example, reaction-diffusion equations, models of convection etc., see
\cite{C.E.}). In particular we will not assume that the attractor is
compact but that it is translation invariant. We will also assume a
space time invariant ergodic measure is given. In other words, our
results apply to $\R \times \R$ actions satisfying the hypothesis
given below. In particular, we assume that the semi-flow on the
function space is a flow when restricted to the attractor of the
system. This follows for example from the analyticity in time of the
solutions of the evolution equations.  The procedure described above
differs with the more standard approach to the space time entropy
which uses boxes of roughly the same size in the space and time
direction. It is however more natural from the point of view of the
definition of the attractors of such systems.

In order to open the possibility of using other type of complexities,
we have tried to isolate the properties we need without reference to a
particular example, although the Kolmogorov complexity satisfies all
the requirements. In order to simplify the proofs, we only discuss the
case of one space dimension, although most results extend easily to
higher dimension.

The paper is organised as follows. In section \ref{results} we
first state the required hypothesis on the dynamical system and on
the complexity, and show that these hypothesis are satisfied by
Kolmogorov complexity. We then state the main results. In section
\ref{proof-1} we prove that under these hypothesis one can define
a complexity per unit time and unit volume. This is done following
the scheme briefly mentioned above of fixing first a finite
precision and removing it only at the end. In section
\ref{proof-2} we prove a variational principle which shows that in
the concrete examples of extended systems studied up to now, the
complexity we have defined is finite. In fact, we show that
computing the supremum of the complexity for functions in the
supports of the invariant measures of the system, one obtains the
topological entropy defined in \cite{C.E.2.}.

\section{Settings and results} \label{results}

Let $\F$ be a set of real functions defined on $\R$ and consider
the following actions on $\F$: the space translation
$$\R \ni y \mapsto (\zeta_y u)(x):= u(x+y)$$
and a flow of time evolution $\varphi_t : \F \to \F$ defined for
$t\in \R$. We assume that the two actions commute.

We assume that the set $\F$ is endowed with a translation
invariant metric $d$ and that there exists a probability measure
$\mu$ on $\F$, such that $\mu$ is invariant and ergodic with
respect to the $(\zeta,\varphi)$ action.

We make the following assumptions on the set $\F$ and the flow
$\varphi$. Let us assume that for any interval $\Lambda \subset
\R$ the set
$$\F|_\Lambda := \set{g:\Lambda \to \R \ :\ \exists \ f\in \F
\hbox{ with } f|_\Lambda \equiv g}$$ is endowed with a metric
$d|_\Lambda$ such that $(\F,d)$ is the projective limit of
$(\F|_\Lambda,d|_\Lambda)$ as $|\Lambda|\to \infty$. If for
example $\F \subset C^0_b(\R)$, the set of real bounded continuous
functions on $\R$, and $d$ is the sup-norm, then for every
$\Lambda \subset \R$ we have $d|_\Lambda(g_1,g_2)=\sup_{x\in
\Lambda} |g_1(x)-g_2(x)|$. We assume that for any interval
$\Lambda \subset \R$ we have
\begin{equation*}
\hbox{{\bf(A1)}} \qquad |\Lambda|< \infty \ \Longrightarrow \
\F|_\Lambda \hbox{ is pre-compact}
\end{equation*}
By assumption \textbf{(A1)}, for any $\eps>0$ and $|\Lambda|<
\infty$ we can define the set
\begin{equation} \label{cylinders}
\C^\eps_\Lambda = \set{\hbox{finite open coverings of
$\F|_\Lambda$ with balls of radius $< \eps$}}
\end{equation}
and we denote by $\U^\eps_\Lambda$ an element of
$\C^\eps_\Lambda$. Fixed two finite intervals $\Lambda_1$ and
$\Lambda_2$ with disjoint interiors let $\Lambda$ be the union
$\Lambda:= \Lambda_1 \cup \Lambda_2$, then we assume that
\begin{itemize}
\item[{\bf (A2)}] there exists an integer $q$ depending only on
the metric $d$ such that, for any $\U^\eps_{\Lambda_1} \in
\C^\eps_{\Lambda_1}$ and $\U^\eps_{\Lambda_2} \in
\C^\eps_{\Lambda_2}$ and two balls $B_1 \in \U^\eps_{\Lambda_1}$
and $B_2 \in \U^\eps_{\Lambda_2}$, either the intersection $B_1
\cap B_2$ is empty or can be covered by $q$ balls of a covering
$\U^\eps_\Lambda \in \C^\eps_\Lambda$.
\end{itemize}

The last assumption on the system is about the separation speed of
two nearby functions with time. We assume that there are constants
$\gamma
>0$, $\Gamma >1$ and $C>0$ such that, for any $|\Lambda|< \infty$
and any $\eps>0$ satisfying $diam(\Lambda)> 2 C \eps^{-1}$ and for
any initial conditions $f_1$ and $f_2$ in $\F$ such that
$d|_\Lambda(f_1,f_2) <\eps$, we have
$$\hbox{{\bf(A3)}} \qquad d|_{\Lambda\setminus \set{d(x,\partial
\Lambda)<C\eps^{-1}(t+1)}}(\varphi_t(f_1),\varphi_t(f_2)) < \Gamma
e^{\gamma t} \eps$$ for any $t\in (0,C^{-1}diam(\Lambda)\eps)$
(cfr. \cite{C.E.2.}).

Under these assumptions a notion of topological entropy for the
flow $\varphi$ has been defined in \cite{C.E.2.}. Let
\begin{equation} \label{maxim-number}
N^\eps_\Lambda(T):= \max \set{\card(S^\eps_\Lambda(T))\ :\
S^\eps_\Lambda(T) \hbox{ is made of
$(\Lambda,T,\eps)$-distinguishable orbits}}
\end{equation}
where we say that $f$ and $g$ in $\F|_\Lambda$ have
$(\Lambda,T,\eps)$-indistinguishable orbits up to time $T$ and
with resolution $\eps$ if
$$d|_\Lambda(\varphi_t(f),\varphi_t(g)) < \eps \qquad \forall \ t\in
(0,T)$$ In \cite{C.E.2.}, \cite{C.E.3.} and \cite{zelik} it is
proved that
\begin{equation} \label{topol-ent}
h_{top} := \lim_{\eps \searrow 0} \ \lim_{|\Lambda|\to \infty} \
\frac{1}{|\Lambda|} \ \lim_{T\to \infty} \ \frac{\lo
N^\eps_\Lambda(T)}{T}
\end{equation}
exists and is finite under some additional assumptions, in fact it
is bounded by $\gamma D_{up}$, where $D_{up}$ is called the upper
local dimension per unit length of the set $\F$ in \cite{C.E.2.}
or capacity per unit length in \cite{kt}.

The aim of this paper is to introduce a measure of the complexity
of the action of the flow $\varphi$ on $\F$ which would be the
analogous of the metric entropy for dynamical systems. To this aim
we need to define a notion of complexity. Our definition is
inspired by the notion of Kolmogorov complexity (\cite{li-vit}).

Let $\A^*$ be the set of finite words on a finite alphabet $\A$,
and for a word $s$ let us denote by $|s|$ its length. We say that
$K:\A^* \to \R^+$, defined for any alphabet $\A$, is a
\emph{\virg{good} complexity function} if it satisfies the
following hypothesis \textbf{(H1)}-\textbf{(H4)}.

The first hypothesis is about the behaviour of the complexity
function on sub-words and a sub-additivity property. Let $s=uv$ be
the concatenation of two words $u$ and $v$, then
$$\hbox{{\bf(H1.a)}} \qquad  K(u) \le K(s) + \lo |u| + const$$
for a constant independent on $s$ and $u$. Moreover let us assume
that there exists a function $h:\N \to \R^+$ satisfying
$\lim_{n\to \infty} \frac{h(n)}{n}=0$ such that
$$\hbox{{\bf(H1.b)}} \qquad  K(s) \le K(u) + K(v) + h(|u|) + h(|v|)$$

Let now $\A_1$ and $\A_2$ be two different alphabets, with $r_i:=
\card(\A_i)$. Moreover let $\tilde \A$ be an alphabet with
$\card(\tilde \A)=q r_1 r_2$ for some integer number $q\ge 1$, and
we assume that there exists a surjective map $\pi: \tilde \A \to
\A_1 \times \A_2$, with coordinate maps $\pi_1$ and $\pi_2$ on
$\A_1$ and $\A_2$, respectively. Let $s \in \tilde \A^*$ and
$\pi_i(s) \in \A_i^*$ be its projections. Then
$$\hbox{{\bf(H2.a)}} \qquad  K(\pi_i(s)) \le K(s) + const$$
$$\hbox{{\bf(H2.b)}} \qquad  K(s) \le K(\pi_1(s)) + K(\pi_2(s)) +
|s| \lo q + const$$ where the constants are independent on $s$.

The third hypothesis is an estimate on $K$ that comes from
observations by Shannon for his definition of information content
(\cite{shannon}). Let $E \subset \A^* \times \N$ be a recursively
enumerable set (for a definition see for example \cite{li-vit}),
and for any $n\in \N$ let $L_n:= \set{s\in \A^*\ :\ (s,n) \in E}$
be a set of finite cardinality. Then we assume that for all $n\in
\N$ it holds
$$\hbox{{\bf(H3)}} \qquad K(s) \le \lo(\card(L_n)) + \lo n +
const \qquad \forall \ s\in L_n$$ where the constant only depends
on the set $E$.

Finally we ask for a relation between the bound of the complexity
on a set of words and the cardinality of this set. We assume that
$$\hbox{{\bf(H4)}} \qquad  \card \set{s\in \A^* \ :\ K(s) < c}
\le 2^c \qquad \forall c\in \R$$

By using a \virg{good} complexity function let us now define the
complexity of the flow $\varphi_t$.

Consider a fixed probability measure $\mu$ which is invariant and
ergodic for the action of $(\zeta,\varphi)$. For a given $\eps >0$
and an interval $\Lambda \subset \R$ with $|\Lambda|< \infty$, we
consider on $\F|_\Lambda$ the set of coverings $\C^\eps_\Lambda$.
We will use such coverings to code the orbit of a function $f\in
\F$ under $\varphi$. To this aim, we introduce a time step
$\tau>0$ and consider the orbits
$(f,\varphi_\tau(f),\varphi_{2\tau}(f),\dots)$. By the method of
symbolic dynamics we can associate to an orbit
$(\varphi_{j\tau}(f))_{j=0}^{n-1}$ a set of $n$-long words
$\psi(f,n,\U^\eps_\Lambda)$ on a finite alphabet
$\A=\A(\U^\eps_\Lambda)=\set{1,\dots,\card(\U^\eps_\Lambda)}$. If
we denote $\U^\eps_\Lambda:=
\set{U_1,\dots,U_{\card(\U^\eps_\Lambda)}}$, we
define\footnote{Since the covering is made by open sets, the
methods of symbolic dynamics give more than one word.}
$$\psi(f,n,\U^\eps_\Lambda):= \set{\omega_0^{n-1} \in
\A(\U^\eps_\Lambda)\ :\ \varphi_{j\tau}(f) \in U_{\omega_j}\
\forall\ j=0,\dots,n-1}$$
In the same way we can define in the
general case $\psi(\varphi_{m\tau}(f),n-m,\U^\eps_\Lambda)$ as the
set of possible symbolic representations of the orbit
$(\varphi_{j\tau}(f))_{j=m}^{n-1}$. At this point we can use a
complexity function $K$ to define
\begin{equation} \label{compl-m-n}
  K(f,\tau,\U^\eps_\Lambda,m,n) := \min \set{K(\omega_m^{n-1})\ :\
  \omega_m^{n-1} \in \psi(\varphi_{m\tau}(f),n-m,\U^\eps_\Lambda)}
\end{equation}
To simplify notations, for $m=0$ we will write
\begin{equation} \label{compl-n}
  K(f,\tau,\U^\eps_\Lambda,n) := \min \set{K(\omega_0^{n-1})\ :\
  \omega_0^{n-1} \in \psi(f,n,\U^\eps_\Lambda)}
\end{equation}
We can define the asymptotic linear rate of increase in $n$ by
\begin{equation} \label{compl-tau}
  K(f,\tau,\U^\eps_\Lambda) := \lim_{n\to \infty}
  \frac{K(f,\tau,\U^\eps_\Lambda,n)}{n}
\end{equation}
To get rid of the dependence on the covering we define
\begin{equation} \label{compl-L}
K(f,\tau,\eps,\Lambda):= \inf \set{K(f,\tau,\U^\eps_\Lambda)\ :
  \U^\eps_\Lambda \in \C^\eps_\Lambda}
\end{equation}
The next step will be to study the asymptotic rate of increase in
$|\Lambda|$. We restrict ourselves to a class of intervals defined
as follows.

\begin{defn} \label{inter-admiss}
A sequence of sets $\La=\set{\Lambda_k}_k$ is called
\emph{admissible} if $\Lambda_k=[a_k,b_k]$ for two sequences
$\set{a_k}_k$ and $\set{b_k}_k$ satisfying $a_k < b_k$ for all
$k\ge 1$ and
\begin{eqnarray}
& \lim\limits_{k\to \infty} (b_k - a_k)  = +\infty
\label{cond-succ-1} \\
& l_a:=\liminf\limits_{k\to \infty}\ \frac{b_k-a_k}{\max
\set{a_k,0}} > 0 \label{cond-succ-2} \\
& l_b:= \liminf\limits_{k\to \infty}\ \frac{b_k-a_k}{-\min
\set{b_k,0}} > 0 \label{cond-succ-3}
\end{eqnarray}
\end{defn}

\noindent Intuitively, this definition says that these sequences
do not move too fast to the left or to the right.

If $\La$ is an admissible sequence of sets, let us define
\begin{equation} \label{compl}
  K_\mu(f,\tau,\eps) := \lim_{k \to \infty}
  \frac{K(f,\tau,\eps,\Lambda_k)}{|\Lambda_k|}
\end{equation}

Given these definitions, we will prove that

\begin{thm} \label{thm-1}
For a given ergodic probability measure $\mu$, if the complexity
function $K$ satisfies \textbf{(H1)} and \textbf{(H2)}, the limits
in (\ref{compl-tau}) and (\ref{compl}) exist almost surely and
$K(f,\tau,\eps)$ is almost surely equal to a constant
$K_\mu(\tau,\eps)$ not depending on the admissible sequence $\La$
of sets. Moreover the function $K_\mu(\tau,\eps)$ is not
decreasing in $\eps$, hence the limit
$$K_\mu(\tau) := \lim_{\eps \to 0} K_\mu (\tau,\eps)$$
exists and moreover there exists a constant $K_\mu$ such that for
all $\tau
>0$
$$\frac{K_\mu(\tau)}{\tau} = K_\mu$$
\end{thm}

\begin{thm} \label{thm-2}
If the complexity function satisfies also \textbf{(H3)} and
\textbf{(H4)}, then
$$\sup \set{K_\mu \ :\ \mu\ \hbox{ invariant probability
measures}} = h_{top}$$ where $h_{top}$ is defined in
(\ref{topol-ent}).
\end{thm}

Before giving the proofs of these theorems, we recall that for a
finite word $s\in \set{0,1}^*$, the Kolmogorov complexity or
Algorithmic Information Content of $s$ is defined as
\begin{equation} \label{aic}
C(s):= \min \set{|w| \ :\ w\in \set{0,1}^*,\ \ U(w)=s}
\end{equation}
where $|\cdot|$ denotes the length of a word, and $U$ is a
universal Turing machine. For more details we refer to
\cite{li-vit}.

\begin{thm} \label{thm-3}
The Kolmogorov complexity satisfies hypotheses
\textbf{(H1)}-\textbf{(H4)}.
\end{thm}

\noindent {\bf Proof.} We recall that the translation of a finite
word from the binary alphabet to any other finite alphabet $\A$
requires only a constant amount of information content not
dependent on the word. Hence we assume that these constants are
included in the hypotheses \textbf{(H1)}-\textbf{(H3)}.

Hypothesis \textbf{(H1)} and \textbf{(H2)} follow from
\cite{li-vit}, equation (2.2) and arguments used in \cite{li-vit},
section 2.1.2.

Hypothesis \textbf{(H3)} is a corollary of Theorem 2.1.3 in
\cite{li-vit}.

Hypothesis \textbf{(H4)} is Theorem 2.2.1 in \cite{li-vit}. \qed

\section{Proof of Theorem \ref{thm-1}} \label{proof-1}

Let us consider any fixed probability measure $\mu$ which is
invariant and ergodic for the action of $(\zeta,\varphi)$.

The first part of the proof relies on the application of arguments
related to the sub-additivity property to define the quantities in
(\ref{compl-tau}) and (\ref{compl}).

Let $X=(X_{m,n})$ be a family of real random variables with
indexes $m,n \in \N$. We recall that $X$ is \emph{almost
subadditive} if there exists a family of random variables
$U=(U_j)$, with $j\in \N$, defined on the same probability space
of $X$ such that
\begin{equation} \label{subaddit}
X_{m,n} \le \sum_{i=1}^{k-1} \ (X_{j_i,j_{i+1}} + U_{j_{i+1}-j_i})
\end{equation}
for all $1\le m < n$ and all partitions $m=j_1< j_2<\dots < j_k
=n$. In \cite{schurger} the following result is proved
\begin{thm}[\cite{schurger}] \label{teo sub}
Let $X$ and $U$ be jointly stationary and let $X$ be almost
subadditive with respect to $U$. Assume that $X^+_{0,1} \in L^1_+$
and that there exists an increasing sequence of integers $(m_k)_k$
with $m_1\ge 1$ such that
\begin{equation} \label{teo sub 1}
\liminf_{k\to \infty}\ \frac{X_{0,n+m_k}}{n+m_k}\ \ge \
\liminf_{k\to \infty}\ \frac{X_{0,m_k}}{m_k} \quad \hbox{ almost
surely}
\end{equation}
for all $n\ge 1$ and
\begin{equation} \label{teo sub 2}
\lim_{k\to \infty}\ \frac{U_{m_k}}{m_k} = 0 \quad \hbox{ almost
surely}
\end{equation}
Then
$$\lim_{k\to \infty}\ \frac{X_{0,m_k}}{m_k} = \bar x \quad \hbox{ exists
almost surely}$$ with $-\infty \le \bar x < \infty$ almost surely.
\end{thm}
We first apply this theorem to $K(f,\tau,\U^\eps_\Lambda,n)$ as
defined in (\ref{compl-n}), identifying $X_{m,n}$ with
$K(f,\tau,\U^\eps_\Lambda,m,n)$. Then, since for all $\omega \in
\A^\N$ it holds $K(\omega_0) \le const$, we have
$$K(f,\tau,\U^\eps_\Lambda,1) = \min \set{K(\omega_0)\ :\ \omega_0
\in \psi(f,1,\U^\eps_\Lambda)} \in L^1$$ Moreover by
\textbf{(H1.a)} we have
$$K(f,\tau,\U^\eps_\Lambda,k) \le
K(f,\tau,\U^\eps_\Lambda,n+k) + \lo k + const$$ for all $n\ge 1$
and all $f\in \F$, hence condition (\ref{teo sub 1}) of the
previous theorem is satisfied with $m_k =k$. We now show the
sub-additivity property with respect to a family of random
variables.

\begin{lem} \label{lem sub 1}
For any fixed $\tau>0$, $\eps>0$, $|\Lambda|<\infty$ and
$\U^\eps_\Lambda \in \C^\eps_\Lambda$, the family
$(K(f,\tau,\U^\eps_\Lambda,m,n))_{m,n}$ is almost subadditive with
respect to the family of functions $h=h(j)$ defined in
\textbf{(H1.b)}.
\end{lem}

\noindent {\bf Proof.} Without loss of generality we can assume
$m=0$ because of stationarity. Let us fix a function $f\in \F$.
For all $\omega_0^{n-1} \in \psi(f,n,\U^\eps_\Lambda)$ by
\textbf{(H1.b)} we have
$$K(\omega_0^{n-1}) \le \sum_{i=1}^{k-1} (K(\omega_{j_i}^{j_{i+1}}) +
h(j_{i+1}-j_i))$$ for any partition $0=j_1< j_2 <\dots < j_k=n-1$.
Fixed any such partition, let $\bar \omega_{j_i}^{j_{i+1}},
i=1,\dots,k-1$, be a collection of finite words such that
$$K(f,\tau,\U^\eps_\Lambda,j_i,j_{i+1})= K(\bar
\omega_{j_i}^{j_{i+1}})$$ Then, if we denote by $\bar
\omega_0^{n-1}$ the concatenation
$$\bar \omega_0^{n-1} := \bar \omega_{j_1}^{j_{2}} \ \bar
\omega_{j_2}^{j_{3}} \dots \bar \omega_{j_{k-1}}^{j_k}$$ it holds
$$\sum_{i=1}^{k-1} (K(\bar \omega_{j_i}^{j_{i+1}}) +
h(j_{i+1}-j_i)) \ge K(\bar \omega_0^{n-1}) \ge
K(f,\tau,\U^\eps_\Lambda,n)$$ hence the sub-additivity property is
proved. \qed

\vskip 0.5cm Since condition (\ref{teo sub 2}) is verified by the
function $h(n)$ and the probability measure $\mu$ is invariant, we
can apply Theorem \ref{teo sub} to obtain that the limit
$K(f,\tau,\U^\eps_\Lambda)$ exists and is finite $\mu$ almost
surely. Let us denote by $Y^\tau_{\U^\eps_\Lambda} \subset \F$ a
set with $\mu((Y^\tau_{\U^\eps_\Lambda})^c)=0$ on which the limit
exists. Then we define
\begin{equation} \label{compl-U}
K(f,\tau,\U^\eps_\Lambda) := \left\{
\begin{array}{ll}
\lim\limits_{n\to \infty} \ \frac{K(f,\tau,\U^\eps_\Lambda,n)}{n}
&
\hbox{ if $f\in Y^\tau_{\U^\eps_\Lambda}$} \\[0.5cm]
+\infty & \hbox{ otherwise}
\end{array} \right.
\end{equation}
We can then prove

\begin{lem} \label{lem sub inf}
There exists a set $Y^\tau_{\eps,\Lambda} \subset \F$ with
$\mu((Y^\tau_{\eps,\Lambda})^c)=0$ such that
$$K(f,\tau,\eps,\Lambda):= \inf \set{K(f,\tau,\U^\eps_\Lambda)\ :
  \U^\eps_\Lambda \in \C^\eps_\Lambda}$$
is well defined and finite for all $f\in Y^\tau_{\eps,\Lambda}$.
Moreover there exists a sequence $\{\tilde \V_s\}_s$ of coverings
in $\C^\eps_\Lambda$ such that
\begin{equation} \label{succ-cov-1}
\lim_{s\to \infty}\ K(f,\tau,\tilde \V_s) = K(f,\tau,\eps,\Lambda)
\qquad \forall \ f\in Y^\tau_{\eps,\Lambda}
\end{equation}
and the sequence $\{K(f,\tau,\tilde \V_s)\}_s$ is non-increasing
for all $f\in Y^\tau_{\eps,\Lambda}$.
\end{lem}

\noindent {\bf Proof.} We first restrict to a countable set of
coverings $\D^\eps_\Lambda \subset \C^\eps_\Lambda$. Let
$G=\set{g_j} \subset \F$ be a countable set dense in
$\F|_\Lambda$, and define $\D^\eps_\Lambda$ as the set of finite
open coverings
\begin{equation} \label{cylinders-count}
\D^\eps_\Lambda := \set{\V^\eps_\Lambda \in \C^\eps_\Lambda \ :\
\hbox{the centers are in $G$ and radii are rational}}
\end{equation}
Restricting to coverings $\V^\eps_\Lambda \in \D^\eps_\Lambda$ we
can define
\begin{equation} \label{compl-L-count}
\tilde K(f,\tau,\eps,\Lambda) := \inf
\set{K(f,\tau,\V^\eps_\Lambda)\ :
  \V^\eps_\Lambda \in \D^\eps_\Lambda}
\end{equation}
on the set $Y^\tau_{\eps,\Lambda} \subset \F$ with
$\mu((Y^\tau_{\eps,\Lambda})^c)=0$ defined by
$$Y^\tau_{\eps,\Lambda} := \bigcap_{\V^\eps_\Lambda \in \D^\eps_\Lambda} \
Y^\tau_{\V^\eps_\Lambda}$$ We now show that $\tilde
K(f,\tau,\eps,\Lambda)$ is equal to $K(f,\tau,\eps,\Lambda)$ on
$Y^\tau_{\eps,\Lambda}$. To this aim it is enough to prove that
for any $\U^\eps_\Lambda \in \C^\eps_\Lambda$ there exists
$\V^\eps_\Lambda \in \D^\eps_\Lambda$ such that
\begin{equation} \label{lemma-part}
K(f,\tau,\V^\eps_\Lambda) \le K(f,\tau,\U^\eps_\Lambda)
\end{equation}
Indeed from this and (\ref{compl-U}), on $Y^\tau_{\eps,\Lambda}$
we have that
$$\tilde K(f,\tau,\eps,\Lambda) \le K(f,\tau,\eps,\Lambda)$$
and the other inequality is obtained by using $\D^\eps_\Lambda
\subset \C^\eps_\Lambda$.

Let us now prove (\ref{lemma-part}). Let $\U^\eps_\Lambda =
\set{U_1,\dots,U_c}$ be a covering in $\C^\eps_\Lambda$ with
$c=\card(\U^\eps_\Lambda)$, and define $r=\max \set{r(U_j)\ :\
j=1,\dots,c}<\eps$ where $r(U_j)$ is the radius of the ball $B_j$.
Then by density of the set of functions $G$ in $\F|_\Lambda$, we
can find a covering $\V^\eps_\Lambda \in \D^\eps_\Lambda$ with
balls $\set{V_j}_{j=1,\dots,c}$ such that $U_j \subset V_j$ for
all $j=1,\dots,c$. Indeed it is enough to choose balls $V_j$ with
centres in functions of the set $G$ at distances less than
$\frac{\eps-r}{2}$ from the centres of the balls $U_j$.

For this choice of coverings for all $n\ge 1$ it holds
$$\psi(f,n,\U^\eps_\Lambda) \subset \psi(f,n,\V^\eps_\Lambda)$$
hence from (\ref{compl-n}) and (\ref{compl-U}) it follows that for
all $f\in Y^\tau_{\eps,\Lambda}$
$$K(f,\tau,\V^\eps_\Lambda)= \lim_{n\to \infty}
\frac{K(f,\tau,\V^\eps_\Lambda,n)}{n} \le \lim_{n\to \infty}
\frac{K(f,\tau,\U^\eps_\Lambda,n)}{n}=K(f,\tau,\U^\eps_\Lambda)$$

We now prove the second part of the assertion. Let us consider an
enumeration of the coverings in $\D^{\eps}_\Lambda=\set{\V_j}_j$,
then we define
$$\tilde \V_s := \bigwedge_{1\le j\le s}\ \V_j$$
where, for two finite open coverings $\U$ and $\V$, by $\U\wedge
\V$ we denote the finite open covering which contains all the
balls of $\U$ and $\V$. By definition, it is clear that $\V_s \in
\D^{\eps}_\Lambda$ for all $s\ge 1$. Moreover, since $\tilde \V_s$
contains all the balls of the coverings $\V_1,\dots,\V_s$, we have
that, modulo a renumbering of the balls of $\tilde \V_s$,
$\psi(f,n,\V_j) \subset \psi(f,n,\tilde \V_s)$ for all
$j=1,\dots,s$ and all $f\in \F$. Hence for all $j=1,\dots,s$
$$K(f,\tau,\tilde \V_s,n) \le K(f,\tau,\V_j,n) +const \qquad \forall
\ f\in \F$$
where the constant is independent on the length $n$ of the
symbolic words. Dividing by $n$ and taking the limit as $n\to
\infty$, we obtain for all $f\in \F$
$$\tilde K(f,\tau,\eps,\Lambda)\le \liminf\limits_{s\to \infty}\
K(f,\tau,\tilde \V_s) \le \limsup\limits_{s\to \infty}\
K(f,\tau,\tilde \V_s) \le \tilde K(f,\tau,\eps,\Lambda)$$ where
the first two inequalities come from the definition of upper and
lower limit. Hence (\ref{succ-cov-1}) is proved.

By the same argument as above, it is immediate to verify that for
all $f\in \F$ the sequence $\{ K(f,\tau,\tilde \V_s) \}_s$ is
non-increasing. Hence the lemma is proved. \qed

\vskip 0.5cm The next step is to show the existence of the limit
in (\ref{compl}) for an admissible sequence of intervals to define
$K(f,\tau,\eps)$. We need the following general lemma

\begin{lem} \label{lemma-succ}
Let $T:(X,\nu) \to (X,\nu)$ be a measure preserving invertible
transformation of a probability space $(X,\nu)$. Let $\vartheta$
and $\xi$ be two real functions on $X$ in the space $L^1(X,\nu)$
and let $\xi(x) \ge 0$ for all $x \in X$. Then there exists a set
$Y\subset X$ with $\nu(Y^c)=0$ such that for any sequences
$\set{a_k}_k$ and $\set{b_k}_k$ of integers satisfying conditions
(\ref{cond-succ-1})-(\ref{cond-succ-3}) we have
\begin{equation} \label{lemma-succ-1}
\bar \vartheta (x) := \lim\limits_{k\to \infty}\ \frac{1}{b_k-a_k}
\sum\limits_{j=a_k}^{b_k-1} \vartheta(T^j(x))
\end{equation}
exists, is finite for all $x\in Y$ and it is in $L^1(X,\nu)$.
Moreover it satisfies
\begin{equation} \label{lemma-succ-3}
\int_X \bar \vartheta (x) d\nu = \int_X \vartheta (x) d\nu
\end{equation}
For the function $\xi$ we have
\begin{equation} \label{lemma-succ-2}
\lim\limits_{k\to \infty}\ \frac{\xi(T^{b_k}(x))}{b_k-a_k} =0
\end{equation} for all $x\in Y$.
\end{lem}

This result is in the spirit of results in \cite{jun-ste} and
\cite{jun-ros}, where it is proved that we cannot ask for weaker
conditions on the sequences $\set{a_k}_k$ and $\set{b_k}_k$.
However we could not relate directly our lemma to their results,
hence in the appendix we give a proof.

We will use this lemma for the space translation action to show
that there exists a set $Y^\tau_\eps\subset \F$, with
$\mu((Y^\tau_\eps)^c)=0$, such that the limit along any admissible
sequence of intervals $\La=\set{\Lambda_k}$
$$\lim_{k\to \infty}\
\frac{K(f,\tau,\eps,\Lambda_k)}{|\Lambda_k|}$$ exists and is
finite for all $f\in Y^\tau_\eps$. The ergodicity of the measure
$\mu$ will imply that this limit is independent on $f$ and it is a
constant $K_\mu(\tau,\eps)$. Moreover from the proof it will
follow that the limit does not depend on the admissible sequence
$\La$ of sets, indeed it will be given by (\ref{ergodico-p}).

We will study separately the superior and the inferior limits. For
the superior limit we use the functions $K(f,\tau,\eps,[0,p])$ for
$p\in \N$. We first prove that we can apply Theorem \ref{teo sub}
to this sequence of functions. We start by verifying that
$K(f,\tau,\eps,[0,1]) \in L^1$. For any $\V^\eps_{[0,1]} \in
\D^\eps_{[0,1]}$ and the associated alphabet $\A$, we can write
$$K(f,\tau,\V^\eps_{[0,1]},n) \le n \ (\max\limits_{\alpha \in \A}
K(\alpha) + const)$$ hence for all $f\in Y_{[0,1]}$
\begin{equation} \label{limitat-in-L}
K(f,\tau,\eps,[0,1]) \le \inf \set{\max\limits_{\alpha \in
\A(\V^\eps_{[0,1]})} K(\alpha)\ :\ \V^\eps_{[0,1]} \in
\D^\eps_{[0,1]}} + const
\end{equation}
Assumption \textbf{(A1)} implies that $K(f,\tau,\eps,[0,1]) \in
L^1$. Note that the bound depends only on the length of the
interval $\Lambda=[0,1]$.

Let $\Lambda_1$ and $\Lambda_2$ be two fixed intervals with
disjoint interiors and denote their union $\Lambda := \Lambda_1
\cup \Lambda_2$, let $\tilde \C^\eps_\Lambda \subset
\C^\eps_\Lambda$ be the set of coverings of $\F|_\Lambda$ built as
in \textbf{(A2)} by two coverings $\U^\eps_{\Lambda_1} \in
\C^\eps_{\Lambda_1}$ and $\U^\eps_{\Lambda_2} \in
\C^\eps_{\Lambda_2}$. For any $\tilde \U^\eps_\Lambda \in \tilde
\C^\eps_\Lambda$ we can write, by using \textbf{(H2.a)}
$$\frac{K(f,\tau,\pi_1(\tilde \U^\eps_{\Lambda}),n)}{n\ |\Lambda_1|} \le
\frac{K(f,\tau,\tilde \U^\eps_{\Lambda},n)}{n\ |\Lambda|} \
\frac{|\Lambda|}{|\Lambda_1|} + \frac{const}{n\ |\Lambda_1|}$$
where $\pi_1(\tilde \U^\eps_{\Lambda})$ denotes the
\virg{projection} of the covering $\tilde \U^\eps_{\Lambda}$ onto
$\C^\eps_{\Lambda_1}$.

Applying this argument with $\Lambda_1=[0,p]$ and
$\Lambda_2=[p,m+p]$ for all $m\ge 1$, and by taking the limit as
$n\to \infty$ and the infimum limit on $|\Lambda_1|=p \to \infty$,
we obtain condition (\ref{teo sub 1}) for any fixed
$|\Lambda_2|=m$.

\begin{lem} \label{lem sub 3}
For any fixed $\eps >0$ and $\tau>0$, the family
$(K(f,\tau,\eps,[0,p]))_{p\in \N}$ is almost subadditive with
respect to the constant function $u(|\Lambda|)\equiv \lo q$, where
$q$ is the constant defined in {\rm \textbf{(A2)}}.
\end{lem}

\noindent {\bf Proof.} Let us consider two disjoint intervals
$\Lambda_1=[a,b]$ and $\Lambda_2=[b,c]$ and the union $\Lambda=
[a,c]$. Let us fix a function $f\in \F$ and let $\tilde
\omega_0^{n-1} \in \psi(f,n,\tilde \U^\eps_\Lambda)$ for a
covering $\tilde \U^\eps_\Lambda \in \tilde \C^\eps_\Lambda$. By
\textbf{(H2.b)} we have
$$K(\tilde \omega_0^{n-1}) \le K(\pi_1(\tilde \omega_0^{n-1})) +
K(\pi_2(\tilde \omega_0^{n-1})) + n \lo q +const$$ since the map
$\pi=(\pi_1,\pi_2)$ is surjective. Then by the same argument as in
Lemma \ref{lem sub 1} we have for all $f\in Y_{\tilde
\U^\eps_\Lambda} \cap Y_{\U^\eps_{\Lambda_1}} \cap
Y_{\U^\eps_{\Lambda_2}}$
$$K(f,\tau,\tilde \U^\eps_\Lambda,n) \le
K(f,\tau,\U^\eps_{\Lambda_1},n) + K(f,\tau,\U^\eps_{\Lambda_2},n)
+ n \lo q + const$$ for $\U^\eps_{\Lambda_i}= \pi_i
(\U^\eps_{\Lambda})$. Then we divide by $n$ and take the limit as
$n\to \infty$. These limits exist as proved above, and we get
$$K(f,\tau,\tilde \U^\eps_\Lambda) \le K(f,\tau,\U^\eps_{\Lambda_1})
+ K(f,\tau,\U^\eps_{\Lambda_2}) + \lo q$$ for all coverings
$\U^\eps_{\Lambda_i} \in \C^\eps_{\Lambda_i}$ and the special
covering $\tilde \U^\eps_{\Lambda} \in \tilde \C^\eps_{\Lambda}$
built from the two. For any fixed $\delta >0$ let us choose two
coverings $\U^\eps_{\Lambda_i}$ satisfying
$$K(f,\tau,\U^\eps_{\Lambda_i}) \le K(f,\tau,\eps,\Lambda_i) + \delta$$
then
$$K(f,\tau,\eps,\Lambda_1) + K(f,\tau,\eps,\Lambda_2) + 2 \delta +\lo q \ge
K(f,\tau,\tilde \U^\eps_\Lambda) \ge K(f,\tau,\eps,\Lambda)$$ and
sub-additivity is proved since it holds for all $\delta >0$. \qed

\vskip 0.5cm Since $u(|\Lambda|)= \lo q$ obviously satisfies
condition (\ref{teo sub 2}), we have that there exist a set
$\tilde Y^\tau_\eps \subset \F$ with $\mu((\tilde
Y^\tau_\eps)^c)=0$ on which $K(f,\tau,\eps,[0,p])$ is defined for
all $p\in \N$, and the limit
\begin{equation} \label{ergodico-p}
\lim_{p\to \infty}\ \frac{K(f,\tau,\eps,[0,p])}{p} =: K_\mu
(f,\tau,\eps)
\end{equation}
exists, is finite almost surely and it is in $L^1(\F,\mu)$.
Moreover the limit holds also in $L^1$ and we denote
\begin{equation} \label{ergodico-p-2}
K_\mu(\tau,\eps) := \int_{\F} \ K_\mu (f,\tau,\eps)\ d\mu
\end{equation}
We remark that if the measure $\mu$ is ergodic then $K_\mu
(f,\tau,\eps)$ is almost surely constant and equal to
$K_\mu(\tau,\eps)$.

Following the notation of Lemma \ref{lemma-succ}, we denote
\begin{equation}
\bar K(f,\tau,\eps,[0,p]) := \lim_{N\to \infty}\ \frac{1}{N}
\sum_{j=0}^{N-1}\ K(\zeta_{jp}f,\tau,\eps,[0,p])
\end{equation}
where it exists. Then we prove the following lemma

\begin{lem} \label{limsup-K-L}
For any fixed $\tau$ and $\eps$, there exists a set $\bar
Y^\tau_\eps$ with $\mu((\bar Y^\tau_\eps)^c)=0$ such that for all
$f\in \bar Y^\tau_\eps$ and any admissible sequence of intervals
$\La=\set{\Lambda_k}$, it holds
$$\limsup\limits_{k\to \infty} \
\frac{K(f,\tau,\eps,\Lambda_k)}{|\Lambda_k|} \le
\liminf\limits_{p\to \infty}\ \frac{\bar
K(f,\tau,\eps,[0,p])}{p}$$ If moreover the measure $\mu$ is
ergodic then
$$\limsup\limits_{k\to \infty} \
\frac{K(f,\tau,\eps,\Lambda_k)}{|\Lambda_k|} \le
K_\mu(\tau,\eps)$$
\end{lem}

\noindent {\bf Proof.} Let us consider an admissible sequence of
intervals with integer boundary points. For a fixed integer $p\in
\N$, we can use the sub-additivity property \textbf{(H2.b)} as in
Lemma \ref{lem sub 3} to show that for all $f\in \tilde
Y^\tau_\eps$ we have, by setting $\lfloor \frac{b_k}{p} \rfloor =:
\tilde b_k$ and $\lfloor \frac{a_k}{p} \rfloor +1 =: \tilde a_k$,
\begin{eqnarray*}
K(f,[a_k,b_k]) \le & \sum_{j=\tilde a_k}^{\tilde b_k-1} \
\left(K(f,[jp,(j+1)p]) +\lo q \right)+ \\
& + K(f,[a_k,\tilde a_k p])+ K(f,[\tilde b_k p,b_k]) + 2\lo q
\end{eqnarray*}
where the dependence on $\tau$ and $\eps$ has been ignored to
simplify notations. First of all, by repeating the same argument
we used to prove (\ref{limitat-in-L}), we can prove that there
exists a constant depending only on $p$, see remark after
(\ref{limitat-in-L}), that is a bound from above for
$K(f,[a_k,\tilde a_k p])$ and $K(f,[\tilde b_k p,b_k])$ for all
$f\in \tilde Y^\tau_\eps$. Moreover we can write
$$K(f,[jp,(j+1)p]) = K(\zeta_{jp}f, [0,p])$$
hence
$$\frac{K(f,[a_k,b_k])}{b_k-a_k} \le \frac{1}{b_k-a_k}
\sum_{j=\tilde a_k}^{\tilde b_k-1}\ \left(K(\zeta_{jp}f,[0,p])
+\lo q\right) + \frac{const}{b_k-a_k}$$ We now apply Lemma
\ref{lemma-succ} to the action of the space translation $\zeta$
and with $K(f,[0,p])$ having the role of the $L^1$ function
$\vartheta$. Let $Y_p \subset \F$ be the full measure set given
for $K(f,[0,p])$ by Lemma \ref{lemma-succ}, then we conclude by
(\ref{lemma-succ-1}) that for all $f\in \bar Y^\tau_\eps := \tilde
Y^\tau_\eps \cap (\cap_p Y_p)$, we have
$$\limsup\limits_{k\to \infty} \frac{K(f,[a_k,b_k])}{b_k-a_k} \le
\frac{\bar K(f,[0,p])}{p} + \frac{\lo q}{p}$$ for all $p\in \N$.
Hence we obtain the first part of the assertion.

The second part follows by first applying Lemma \ref{lemma-succ}
to an ergodic measure, from which we get that for all $p\in \N$
$$\bar K(f,[0,p]) = \int_X \ K(f,[0,p]) \ d\mu$$
almost surely. Then we use (\ref{ergodico-p}) and
(\ref{ergodico-p-2}) to conclude.

The result for the sequences $\set{a_k}_k$ and $\set{b_k}_k$
follows by writing
\begin{eqnarray*}
K(f,[a_k,b_k]) & \le K(f,[a_k,\lfloor a_k \rfloor +1]) +
K(f,[\lfloor a_k \rfloor +1,\lfloor b_k \rfloor])+ \\ & +
K(f,[\lfloor b_k \rfloor,b_k]) + 3\lo q
\end{eqnarray*}
and reducing to the above argument. \qed

\vskip 0.5cm To prove a similar result for the inferior limit we
use the following result proved in \cite{derr}

\begin{thm}[\cite{derr}] \label{teo derr}
Let $T:(X,\nu) \to (X,\nu)$ be a measure preserving invertible
transformation of a probability space $(X,\nu)$. Let
$\set{\beta_n}_n$ be a sequence of integrable real functions on
$X$ such that
\begin{equation} \label{derr1}
\inf_n \ \frac{1}{n} \ \int_X\ \beta_n(x)\ d\nu(x) > -\infty
\end{equation}
and for all $n,k$
\begin{equation} \label{derr2}
\beta_{n+k}(x) - \beta_n(x) - \beta_k(T^n(x)) \le h_k(T^n(x))
\end{equation}
for a sequence of functions $\set{h_k}_k$ satisfying $h_k \ge 0$
and $\int_X h_k d\nu \le const$. Then there exists a function
$\vartheta \in L^1(X,\nu)$ such that
\begin{equation} \label{derr3}
\int_X \ \vartheta(x) \ d\nu = \lim_{n\to \infty} \frac{1}{n} \
\int_X\ \beta_n(x)\ d\nu
\end{equation}
and a function $\xi \in L^1(X,\nu)$ such that $\xi \ge 0$ and
\begin{equation} \label{derr4}
\sum_{j=0}^{n-1}\ \vartheta(T^j(x)) \le \beta_n(x) + \xi(T^n(x))
\end{equation}
for almost all $x\in X$ and all $n\in \N$.
\end{thm}

\begin{lem} \label{liminf-K-L}
Under the hypothesis of ergodicity for the probability measure
$\mu$ on $\F$, for any fixed $\tau$ and $\eps$, there exists a set
$\underline Y^\tau_\eps$ with $\mu((\underline Y^\tau_\eps)^c)=0$
such that for all $f\in \underline Y^\tau_\eps$ and any admissible
sequence of intervals $\La=\set{\Lambda_k}$, it holds
$$\liminf\limits_{k\to \infty} \
\frac{K(f,\tau,\eps,\Lambda_k)}{|\Lambda_k|} \ge
K_\mu(\tau,\eps)$$
\end{lem}

\noindent {\bf Proof.} Let us consider first the case of sequences
$\set{a_k}$ and $\set{b_k}$ of integers. We apply Theorem \ref{teo
derr} to the sequence $\beta_p(f):=K(f,\tau,\eps,[0,p])$ for $f\in
\tilde Y^\tau_\eps$ and $p\in \N$. Indeed condition (\ref{derr1})
is easily verified since $K(f,[0,p])$ is non-negative\footnote{To
simplify notations we neglect the dependence on $\tau$ and
$\eps$.}, and condition (\ref{derr2}) is the subadditive property
we proved in Lemma \ref{lem sub 3}. Hence there exists a full
measure set $Y(\xi,\vartheta) \subset \F$ on which (\ref{derr4})
is verified for two given functions $\vartheta$ and $\xi$.
Moreover from (\ref{ergodico-p}) and (\ref{ergodico-p-2}), for an
ergodic measure $\mu$ we have that
$$\int_{\F} \ \vartheta(f) \ d\mu = K_\mu(\tau,\eps)$$
Now let us define the set
$$Y_\zeta(\xi,\vartheta) := \bigcap_{n\in \Z} \set{\zeta_n f \ :\
f\in Y(\xi,\vartheta)}$$ since $\mu$ is $\zeta$-invariant it holds
$\mu((Y_\zeta(\xi,\vartheta))^c)=0$. On $Y_\zeta(\xi,\vartheta)$
we can write
$$K(f,[a_k,b_k]) = K(\zeta_{a_k}f,[0,b_k-a_k])$$
and, using (\ref{derr4}),
$$\sum_{j=0}^{b_k-a_k-1}\ \vartheta(\zeta_{j+a_k}(f)) \le
K(f,[a_k,b_k]) + \xi(\zeta_{b_k}(f))$$ At this point we apply
Lemma \ref{lemma-succ} to $\vartheta$ and $\xi$, and we obtain
using the ergodicity of the measure $\mu$
$$K_\mu(\tau,\eps) = \lim_{k\to \infty} \ \frac{1}{b_k-a_k}
\sum_{j=a_k}^{b_k-1}\ \vartheta(\zeta_j(f)) \le
\liminf\limits_{k\to \infty} \frac{K(f,[a_k,b_k])}{b_k-a_k}$$ for
almost all $f\in \F$. Let us call this full measure set
$\underline Y^\tau_\eps$.

As in Lemma \ref{limsup-K-L}, the general case for real sequences
$\set{a_k}$ and $\set{b_k}$ follows by
\begin{eqnarray*}
K(f,[a_k,b_k]) & \ge K(f,[\lfloor a_k \rfloor,\lfloor b_k
\rfloor+1]) - K(f,[\lfloor a_k \rfloor ,a_k]) - \\
& + K(f,[b_k,\lfloor b_k \rfloor+1]) - 3\lo q
\end{eqnarray*}
The lemma is proved. \qed

\vskip 0.5cm Putting together Lemmas \ref{limsup-K-L} and
\ref{liminf-K-L} we have that, for any fixed $\tau$ and $\eps$,
there exists a set $Y^\tau_\eps := \bar Y^\tau_\eps \cap
\underline Y^\tau_\eps \subset \F$ with $\mu((Y^\tau_\eps)^c)=0$,
such that for all $f\in Y^\tau_\eps$ the limit
$$\lim_{k\to \infty} \
\frac{K(f,\tau,\eps,\Lambda_k)}{|\Lambda_k|}= K_\mu(f,\tau,\eps)$$
exists and is finite for all admissible sequence of intervals
$\La=\set{\Lambda_k}$. Moreover, for an ergodic measure $\mu$, for
all $f\in Y^\tau_\eps$ and any admissible sequence of intervals
$\La$, it is equal to the constant $K_\mu(\tau,\eps)$ defined in
(\ref{ergodico-p-2}). Hence the limit in (\ref{compl}) exists.

To finish the proof of Theorem \ref{thm-1}, we first have to prove
that $K_\mu(\tau,\eps)$ is non decreasing in $\eps$. Let $\eps_1 <
\eps_2$, then according to the above arguments, we can define
$K(f,\tau,\eps_1,\Lambda)$ and $K(f,\tau,\eps_2,\Lambda)$ for all
$|\Lambda|<\infty$ as in (\ref{compl-L}) and all $f\in
Y^\tau_{\eps_1} \cap Y^\tau_{\eps_2}$ . Moreover, since
$\C^{\eps_1}_\Lambda \subset \C^{\eps_2}_\Lambda$ we have
$Y^\tau_{\eps_1,\Lambda} \subset Y^\tau_{\eps_2,\Lambda}$ and
$$K(f,\tau,\eps_2,\Lambda) \le K(f,\tau,\eps_1,\Lambda)$$
Restricting to $\Lambda=[a_k,b_k]$ with $a_k < b_k$ for all $k\ge
1$ and satisfying (\ref{cond-succ-1})-(\ref{cond-succ-3}),
dividing by $|\Lambda|=(b_k-a_k)$ and taking the limit as in
(\ref{compl}) gives
$$K_\mu(\tau,\eps_2) \le K_\mu(\tau,\eps_1) < \infty$$ on the set
$Y^\tau_{\eps_1}\subset Y^\tau_{\eps_2}$. Let us now consider a
monotonically vanishing sequence $\set{\eps_k}_k \subset \Q^+$ and
define
$$Y^\tau := \bigcap\limits_{k \ge 1} \ Y^\tau_{\eps_k}$$
Then $\mu((Y^\tau)^c)=0$ and $K_\mu(\tau,\eps)$ is finite on
$Y^\tau$ for all $\eps>0$ and non-decreasing on $\eps$. Hence we
can define $$K_\mu(\tau) = \lim_{\eps \to 0} K_\mu(\tau,\eps)$$ on
the full measure set of functions $Y^\tau$. In Theorem \ref{thm-2}
we prove that it is finite for complexity functions satisfying
also \textbf{(H3)} and \textbf{(H4)}.

Now it remains to prove that $\frac{K_\mu(\tau)}{\tau}$ does not
depend on $\tau$. We will use assumption \textbf{(A3)}.

\begin{lem} \label{lem sub 2}
The full measure set $Y:=Y^\tau$ does not depend on $\tau$, and
there exists a constant $K_\mu$ such that
$\frac{K_\mu(\tau)}{\tau}=K_\mu$ on $Y$ for all $\tau
>0$.
\end{lem}

\noindent {\bf Proof.} Let us fix constants $\gamma,\Gamma,C$ as
in \textbf{(A3)}, an interval $|\Lambda|<\infty$, $\eps >0$ and a
time step $\tau \in (0,C^{-1}diam(\Lambda)\eps)$. Then for all
$\tau'>0$ we denote
$$\eps':=\Gamma e^{\gamma \tau} \eps$$
$$\Lambda':=\Lambda \setminus \set{d(x,\partial
\Lambda)<C\eps^{-1}(\tau+1)}$$ Moreover, given a covering
$\U^\eps_\Lambda \in \C^\eps_\Lambda$, we consider the covering
$\bar \U^{\eps'}_{\Lambda'} \in \C^{\eps'}_{\Lambda'}$ which has
balls with the same centres as those in $\U^{\eps}_{\Lambda}$ and
radius increased by a factor $\Gamma e^{\gamma \tau}$. By
assumption \textbf{(A3)}, for any function $f\in \F$, we can build
a symbolic orbit $\bar \omega_0^{n'-1} \in \psi(f,n',\bar
\U^{\eps'}_{\Lambda'})$ by using the information contained in a
symbolic orbit $\omega_0^{n-1} \in \psi(f,n,\U^\eps_\Lambda)$,
with $n = n' \frac{\tau'}{\tau}$ by defining
$$\bar \omega_{j'} = \omega_{j}\ \hbox{ with } j=
\left\{
\begin{array}{ll}
j'+ \lfloor j' \left( \frac{\tau'}{\tau} -1 \right) \rfloor
& \hbox{ if $\tau' > \tau$} \\[0.3cm]
\lfloor j' \frac{\tau'}{\tau} \rfloor & \hbox{ if $\tau' < \tau$}
\end{array}
\right.
$$
hence
\begin{equation} \label{stima-tau}
K(\bar \omega_0^{n'-1}) \le K(\omega_0^{n-1}) + n'
\frac{\tau'}{\tau} + const
\end{equation}
with a constant dependent only on the complexity function, and the
term $n' \frac{\tau'}{\tau}$ that contains the information we need
each time that the difference between $j$ and $j'$ increases of
one unit.

Let now $\bar \omega_0^{n-1}$ be the symbolic orbit on which the
minimum for the function $K(f,\tau,\U^\eps_\Lambda,n)$ is
attained. Then by (\ref{stima-tau}) we have
$$K(f,\tau',\bar \U^{\eps'}_{\Lambda'},n') \le K(\bar
\omega_0^{n'-1})\le K(f,\tau,\U^\eps_\Lambda,n) + n'
\frac{\tau'}{\tau} + const$$ and dividing by $n$ and taking the
limit for $n\to \infty$ we have
\begin{equation} \label{stima-tau-L}
\frac{K(f,\tau',\bar \U^{\eps'}_{\Lambda'})}{\tau'} \le
\frac{K(f,\tau,\U^{\eps}_{\Lambda})}{\tau} + \frac 1 \tau
\end{equation}
for all $f\in Y^\tau \cap Y^{\tau'}$. At this point, let us fix a
$\delta >0$ and let $\U^{\eps}_{\Lambda}$ be a covering such that
$$K(f,\tau,\U^{\eps}_{\Lambda}) < K(f,\tau,\eps,\Lambda) +\delta$$
then using the induced covering $\bar \U^{\eps'}_{\Lambda'}$ and
(\ref{stima-tau-L}) we have
\begin{equation} \label{stima-tau-e}
\frac{K(f,\tau',\eps',\Lambda')}{\tau' |\Lambda|} \le
\frac{K(f,\tau',\bar \U^{\eps'}_{\Lambda'})}{\tau' |\Lambda|} <
\frac{K(f,\tau,\eps,\Lambda) +\delta+1}{\tau |\Lambda|}
\end{equation}
for all $\Lambda=[a_k,b_k]$ as in
(\ref{cond-succ-1})-(\ref{cond-succ-3}). Then the limit for
$|\Lambda| \to \infty$ gives
$$\frac{K_\mu(\tau',\eps')}{\tau'} < \frac{K_\mu(\tau,\eps)}{\tau}
<\infty$$ on $Y^\tau \cap Y^{\tau'}$, where we have used
$\frac{|\Lambda'|}{|\Lambda|}\to 1$ as $|\Lambda|\to \infty$, and
we have suppressed the dependence on the function $f$ because of
the ergodicity of the measure $\mu$.

The final step is the limit for $\eps$, and since $\eps'\to 0$ as
$\eps \to 0$ we have
$$\frac{K_\mu(\tau')}{\tau'} < \frac{K_\mu(\tau)}{\tau}< \infty$$
on $Y^\tau \cap Y^{\tau'}=Y^\tau$.

Repeating the argument interchanging the roles of $\tau$ and
$\tau'$, the lemma is proved. \qed

\vskip 0.5cm Hence Theorem \ref{thm-1} is proved.

\section{Proof of Theorem \ref{thm-2}} \label{proof-2}

Since we proved that $\frac{K_\mu(\tau)}{\tau}=K_\mu$ for all
$\tau>0$, in this proof we can fix $\tau=1$ for simplicity of
notation and drop it from formulas.

The first inequality
\begin{equation} \label{prima-diseg}
\sup \set{K_\mu \ :\ \mu\ \hbox{ invariant probability measures}}
\le h_{top}
\end{equation}
follows by showing that for any $(\zeta,\varphi)$-invariant
probability measure $\mu$, and for any fixed $\eps>0$ and any
finite interval $\Lambda$ it holds
\begin{equation} \label{prima-diseg-scopo}
\int_{\F}\ K(f,\eps,\Lambda) \ d\mu \le \lim_{T\to \infty}
\frac{\lo(N^{\eps/4}_\Lambda(T))}{T}
\end{equation}
where the right hand side is defined as in (\ref{topol-ent}).
Indeed (\ref{prima-diseg-scopo}) implies (\ref{prima-diseg}) just
dividing by $|\Lambda|$, and using the $L^1$ convergence proved in
Theorem \ref{thm-1} for the limits as $|\Lambda|\to \infty$ and
$\eps\to 0$.

To prove (\ref{prima-diseg-scopo}), since for any $f\in \F$ it
holds $K(f,\eps,\Lambda) \le K(f,\U^\eps_\Lambda)$ for any
covering $\U^\eps_\Lambda \in \C^\eps_\Lambda$, it is enough to
prove the following lemma

\begin{lem} \label{lebesgue-top-ent}
There exists a covering $\bar \U^\eps_\Lambda \in \C^\eps_\Lambda$
such that
\begin{equation} \label{scopo-lebes}
\int_{\F}\ K(f,\bar \U^\eps_\Lambda) \ d\mu \le \lim_{n\to \infty}
\frac{\lo(N^{\eps/4}_\Lambda(n))}{n}
\end{equation}
\end{lem}

\noindent {\bf Proof.} Let us consider a covering $\U^\eps_\Lambda
\in \C^\eps_\Lambda$ with balls of radius $\frac \eps 3$, and a
new covering $\bar \U^\eps_\Lambda$ with balls with the same
centres as before and radius $\frac{2\eps}{3}$.

We recall that the Lebesgue number lemma states that for a finite
open covering of a compact metric space, there is a finite number
$\gamma>0$ such that the $\gamma$-neighbourhood of any point is
contained in at least one open set of the covering. The number
$\gamma$ is called the Lebesgue number of the covering. Since
$\F|_\Lambda$ is a metric compact set, the covering $\bar
\U^\eps_\Lambda$ has a finite Lebesgue number, and by its
construction we conclude that its Lebesgue number is not less than
$\frac \eps 4$.

The idea of the proof is to use \textbf{(H3)}, hence we need to
count all possible symbolic orbits. For all $n$ let us consider a
minimal $(n,\frac \eps 4)$-spanning set $\Sigma(n,\frac \eps 4)$
for $\F|_\Lambda$, that is for any $f\in \F$ there exists $g\in
\Sigma(n,\frac \eps 4)$ such that
$$d|_\Lambda (\varphi_k(f),\varphi_k(g))\le \frac \eps 4 \qquad
\forall\ k=0,\dots,n-1$$ For any $g\in \Sigma(n,\frac \eps 4)$ let
us denote by $R_g\subset \F$ the set of functions $f\in \F$ which
are $\frac \eps 4$-spanned by $g$. We can make $\set{R_g}_g$ a
partition of $\F$ just by choosing for each $f\in \F$ only one
function spanning it.

For any function $g\in \Sigma(n,\frac \eps 4)$ we can consider the
sequence of balls $\set{B(\varphi_k(g),\frac \eps 4)}_k$ with
$k=0,\dots,n-1$. At the same time, by our result on the Lebesgue
number of $\bar \U^\eps_\Lambda$, we can associate to each such
sequence of balls a symbolic orbit $\omega_0^{n-1}(g) \in
\psi(g,n,\bar \U^\eps_\Lambda)$. Then we have
\begin{equation} \label{spann-comp}
\int_{\F}\ \frac{K(f,\bar \U^\eps_\Lambda,n)}{n} \ d\mu \le \frac
1 n \ \sum_{g\in \Sigma(n,\frac \eps 4)}\ K(\omega_0^{n-1}(g))\
\mu(R_g)
\end{equation}
Let us now consider for the alphabet $\A=\set{1,\dots,\card(\bar
\U^\eps_\Lambda)}$, the set
$$E:= \bigcup_{n\in \N}\ \bigcup_{g\in \Sigma(n,\frac \eps 4)} \
(\omega_0^{n-1}(g),n) \subset \A^* \times \N$$ It is a recursively
enumerable set, hence we can apply hypothesis \textbf{(H3)}
getting
\begin{equation} \label{applic-h3}
K(\omega_0^{n-1}(g)) \le \lo\left(\sigma\left(n,\frac \eps
4\right)\right)+ \lo n +const
\end{equation}
where $\sigma\left(n,\frac \eps 4\right) := \card
(\Sigma\left(n,\frac \eps 4\right))$. Since this estimate is
uniform on $g$, applying it to (\ref{spann-comp}) we get
\begin{equation} \label{spann-comp-2}
\int_{\F}\ \frac{K(f,\bar \U^\eps_\Lambda,n)}{n} \ d\mu \le
\frac{\left(\lo\left(\sigma\left(n,\frac \eps 4\right)\right)+ \lo
n +const\right)}{n}
\end{equation}
since $\sum_g \mu(R_g)=1$.

To finish the proof of the lemma we use the inequality
$$\sigma\left(n,\frac \eps 4\right) \le N^{\eps/4}_\Lambda(n)$$
which is well known in ergodic theory, see for example \cite{dgs}.
\qed

\vskip 0.5cm Let us now prove the other inequality. For any fixed
$\eps>0$, from the definition of topological entropy
(\ref{topol-ent}) we define\footnote{We recall that we consider
fixed $\tau=1$.}
\begin{equation} \label{h-lambda-eps}
h_\Lambda(\eps):= \lim_{n\to \infty} \ \frac{\lo
N^\eps_\Lambda(n)}{n}
\end{equation}
\begin{equation} \label{h-eps}
h(\eps):=\lim_{|\Lambda| \to \infty} \
\frac{h_{\Lambda}(\eps)}{|\Lambda|}
\end{equation}
The quantity $h(\eps)$ is non-decreasing in $\eps$ and its limit
as $\eps \to 0$ is $h_{top}$. With respect to the quantities
defined above, given any fixed $\delta
>0$ there exist $\eps_0>0$, $\lambda_0(\eps)>0$ and
$n_0(\eps,\Lambda)>0$ such that
\begin{equation} \label{h-eps-0}
h_{top} - \delta < h(\eps) < h_{top}  \qquad \ \forall \ \eps <
\eps_0
\end{equation}
\begin{equation} \label{h-eps-lambda}
(h(\eps)-\delta)|\Lambda| < h_\Lambda(\eps) <
(h(\eps)+\delta)|\Lambda|  \qquad \ \forall \ |\Lambda| >
\lambda_0(\eps)
\end{equation}
\begin{equation} \label{h-eps-lambda-n}
2^{n(h_\Lambda(\eps)-\delta|\Lambda|)} < N^\eps_\Lambda(n) <
2^{n(h_\Lambda(\eps)+\delta|\Lambda|)}  \qquad \ \forall \ n >
n_0(\eps,\Lambda)
\end{equation}

We first state a lemma we need in the following

\begin{lem} \label{convergenza-misure}
Let us consider a fixed $\eps>0$ and a finite interval $\Lambda$.
If $\set{\rho_j}_j$ is a sequence of probability measure on $\F$,
invariant with respect to the time evolution $\varphi_1$, then
there exists a sub-sequence $\set{j_h}_h$ such that
$\set{\rho_{j_h}}_h$ is weakly convergent to a probability measure
$\rho$, invariant with respect to $\varphi_1$, and
$$\limsup\limits_{h\to \infty} \ \int_{\F}\ K(f,\eps,\Lambda) \
d\rho_{j_h} \le \ \int_{\F}\ K(f,\eps,\Lambda) \ d\rho$$
\end{lem}

\noindent {\bf Proof.} The existence of the $\varphi_1$-invariant
probability measure $\rho$ follows by the compactness of the space
$\F|_\Lambda$. For simplicity of notations, let us assume that
$\set{\rho_j}_j$ is weakly convergent to $\rho$.

By the monotonicity of the sequence $\{K(f,\tilde \V_s)\}_s$
proved in Lemma \ref{lem sub inf}, we have
\begin{equation} \label{limiti-misure-1}
\int_{\F}\ K(f,\eps,\Lambda) \ d\rho_j = \lim_{s\to \infty}
\int_{\F}\ K(f,\tilde \V_s) \ d\rho_j
\end{equation}
for all $j\in \N$ and also for the $\varphi_1$-invariant measure
$\rho$. By the sub-additive ergodic theorem in $L^1$ it holds
\begin{equation} \label{limiti-misure-2}
\int_{\F}\ K(f,\tilde \V_s) \ d\rho_j = \lim_{n\to \infty}
\int_{\F}\ \frac{K(f,\tilde \V_s,n)}{n} \ d\rho_j = \inf_{n\to
\infty} \int_{\F}\ \frac{K(f,\tilde \V_s,n)}{n} \ d\rho_j
\end{equation}
for all $s\in \N$. Moreover, for each fixed $n\in \N$, the
function $f\mapsto K(f,\tilde \V_s,n)$ is upper semi-continuous,
hence using weak convergence of $\rho_j$ we get for all $s\in \N$
(see for example \cite{billingsley})
\begin{equation} \label{limiti-misure-3}
\inf_{n\to \infty} \int_{\F}\ \frac{K(f,\tilde \V_s,n)}{n} \
d\rho_j \le \inf_{n\to \infty} \int_{\F}\ \frac{K(f,\tilde
\V_s,n)}{n} \ d\rho = \int_{\F}\ K(f,\tilde \V_s) \ d\rho
\end{equation}
where for the last equality we have used (\ref{limiti-misure-2})
for the $\varphi_1$-invariant measure $\rho$. The assertion
follows by putting together (\ref{limiti-misure-2}) and
(\ref{limiti-misure-3}), and by applying (\ref{limiti-misure-1})
to both sides. \qed

\begin{lem} \label{misura-locale}
Given any $\delta>0$, there exists $\eps_0>0$ such that for any
$\eps<\eps_0$ there exists $\lambda_0(\eps)>0$ such that for any
finite interval $|\Lambda|> \lambda_0(\eps)$ it holds: there
exists a probability measure $\mu^\eps_\Lambda$ on $\F$, invariant
with respect to the time evolution $\varphi_1$, such that
$$\frac{1}{|\Lambda|}\ \int_{\F}\ K(f,\eps/2,\Lambda) \ d
\mu^\eps_\Lambda \ge h_{top} - 4\delta$$
\end{lem}

\noindent {\bf Proof.} For any fixed $\delta>0$ let us consider
$\eps_0>0$ as defined for (\ref{h-eps-0}), and let us fix
$\eps<\eps_0$. Referring to (\ref{maxim-number}), for any finite
interval $\Lambda$ let us denote by
$S^\eps_\Lambda(n):=\set{f_j}_{j=1}^{N^\eps_\Lambda(n)}$ the
functions of a maximal set of $(\Lambda,n,\eps)$-different orbits.
On this set we consider the sequence of probability measures on
$\F$ given by
\begin{equation} \label{misure-conc}
\nu^\eps_{\Lambda,n} := \frac{1}{N^\eps_\Lambda(n)} \
\sum_{j=1}^{N^\eps_\Lambda(n)} \ \delta_{f_j}
\end{equation}
where $\delta_\cdot$ denotes the usual Dirac mass. Hence by
definition of the set $S^\eps_\Lambda(n)$, for any open covering
$\U^{\eps/2}_\Lambda \in \C^{\eps/2}_\Lambda$ we have
\begin{equation} \label{integ-conc}
\frac{1}{|\Lambda|}\ \int_{\F}\
\frac{K(f,\U^{\eps/2}_\Lambda,n)}{n} \ d\nu^\eps_{\Lambda,n} =
\frac{1}{|\Lambda|}\ \frac{1}{N^\eps_\Lambda(n)}\
\sum_{j=1}^{N^\eps_\Lambda(n)}\
\frac{K(f_j,\U^{\eps/2}_\Lambda,n)}{n}
\end{equation}
We now use hypothesis \textbf{(H4)}. For any given $\delta >0$ we
have
\begin{equation} \label{applic-h4}
\card \set{f\in S^\eps_\Lambda(n) \ :\ K(f,\U^{\eps/2}_\Lambda,n)
< n |\Lambda| (h(\eps)-3\delta)} \le 2^{n |\Lambda|
(h(\eps)-3\delta)}
\end{equation}
Putting together (\ref{integ-conc}) and (\ref{applic-h4}) we
obtain
\begin{equation} \label{integ-conc-stima}
\frac{1}{|\Lambda|}\ \int_{\F}\
\frac{K(f,\U^{\eps/2}_\Lambda,n)}{n} \ d\nu^\eps_{\Lambda,n} \ge
\frac{\left(N^\eps_\Lambda(n)-2^{n |\Lambda|
(h(\eps)-3\delta)}\right) (h(\eps)-3\delta)}{N^\eps_\Lambda(n)}
\end{equation}

For any fixed $m\in \N$ let us write $n=tm+r$ with $0\le r<m$.
Using the sub-additivity property for the family of functions
$(K(f,\U^{\eps/2}_\Lambda,m,n))_{m,n}$ proved in Lemma \ref{lem
sub 1}, we write for any $f\in \F$
\begin{equation}
\begin{array}{ll}
K(f,\U^{\eps/2}_\Lambda,n) \le & \sum_{j=0}^{t-2} \
K(\varphi_l(f),\U^{\eps/2}_\Lambda,jm,(j+1)m) + \\[0.3cm]
& + K(f,\U^{\eps/2}_\Lambda,l) +
K(f,\U^{\eps/2}_\Lambda,(t-1)m+l,n)+ th(m)
\end{array}
\end{equation}
for all $l=0,\dots,m-1$, where $h(\cdot)$ is the function defined
in \textbf{(H1.b)}. Moreover for all $j=0,\dots,t-2$ and all
$l=0,\dots,m-1$, it holds
$$\int_{\F}\ K(\varphi_l(f),\U^{\eps/2}_\Lambda,jm,(j+1)m) \
d\nu^\eps_{\Lambda,n} = \int_{\F}\ K(f,\U^{\eps/2}_\Lambda,m) \
d(\varphi_{l+jm}^* \nu^\eps_{\Lambda,n})$$ and for all $f\in \F$
the uniform estimate
$$K(f,\U^{\eps/2}_\Lambda,l) +
K(f,\U^{\eps/2}_\Lambda,(t-1)m+l,n) \le 2m \lo
(\card(\U^{\eps/2}_\Lambda)) + const$$ holds. Hence letting
$$\mu^\eps_{\Lambda,n}:= \frac 1 m \ \sum_{l=0}^{m-1} \ \left(
\frac{1}{t-1} \ \sum_{j=0}^{t-2} \ \varphi_{l+jm}^*
\nu^\eps_{\Lambda,n} \right) = \frac{1}{(t-1)m} \
\sum_{i=0}^{(t-1)m-1} \ \varphi_i^* \nu^\eps_{\Lambda,n}$$ we
obtain
\begin{equation} \label{utile-m}
\int_{\F}\ \frac{K(f,\U^{\eps/2}_\Lambda,n)}{n} \
d\nu^\eps_{\Lambda,n} \le \ \frac{(t-1)m}{tm+r}\ \int_{\F}\
\frac{K(f,\U^{\eps/2}_\Lambda,m)}{m} \ d\mu^\eps_{\Lambda,n} +
\frac{th(m)+o(n)}{tm+r}
\end{equation}
Let $\mu^\eps_{\Lambda}$ be an accumulation point for the sequence
of probability measures $\{\mu^\eps_{\Lambda,n}\}_n$ given by
Lemma \ref{convergenza-misure}. It follows that
$\mu^\eps_{\Lambda}$ is a probability measure on $\F$ which is
invariant for the time action $\varphi_1$. Using (\ref{utile-m})
and the upper semi-continuity of the function $f\mapsto
K(f,\U^{\eps/2}_\Lambda,m)$ for all $m\in \N$, we have
\begin{equation} \label{per-tutti-m}
\limsup\limits_{n\to \infty}\ \int_{\F}\
\frac{K(f,\U^{\eps/2}_\Lambda,n)}{n} \ d\nu^\eps_{\Lambda,n} \le
\int_{\F}\ \frac{K(f,\U^{\eps/2}_\Lambda,m)}{m} \
d\mu^\eps_{\Lambda} + \frac{h(m)}{m}
\end{equation}
for all $m\in \N$.

Let now $\lambda_0(\eps)>0$ and $n_0(\eps,\Lambda)>0$ be defined
as in (\ref{h-eps-lambda}) and (\ref{h-eps-lambda-n}). Then from
(\ref{integ-conc-stima}) we have that if $|\Lambda|> \lambda_0$
and $n>n_0$ it holds
$$\frac{1}{|\Lambda|}\ \int_{\F}\
\frac{K(f,\U^{\eps/2}_\Lambda,n)}{n} \ d\nu^\eps_{\Lambda,n} \ge
\left(1 - 2^{-n|\Lambda|\delta} \right) (h(\eps)-3\delta)$$ Hence
for all $m\in \N$, if $\eps< \eps_0$ and $|\Lambda|>
\lambda_0(\eps)$ we have
\begin{equation} \label{per-tutti-m-fine}
\int_{\F}\ \frac{K(f,\U^{\eps/2}_\Lambda,m)}{m} \
d\mu^\eps_{\Lambda} + \frac{h(m)}{m} \ge (h_{top} - 4\delta)
|\Lambda|
\end{equation}
for all coverings $\U^{\eps/2}_\Lambda \in \C^{\eps/2}_\Lambda$,
where we have used (\ref{h-eps-0}). By the sub-additive ergodic
theorem in $L^1$ we have that
$$\lim_{m\to \infty} \ \int_{\F}\ \frac{K(f,\U^{\eps/2}_\Lambda,m)}{m} \
d\mu^\eps_\Lambda = \int_{\F}\ K(f,\U^{\eps/2}_\Lambda) \
d\mu^\eps_\Lambda$$ where $K(f,\U^{\eps/2}_\Lambda)$ is defined as
in (\ref{compl-U}). Hence, since $h(m)=o(m)$, from
(\ref{per-tutti-m-fine}) we have that given any $\delta>0$ there
exists $\eps_0>0$ such that for any $\eps<\eps_0$ there exists
$\lambda_0(\eps)>0$ such that for any $|\Lambda|>\lambda_0$ it
holds
\begin{equation} \label{per-U}
\int_{\F}\ K(f,\U^{\eps/2}_\Lambda) \ d\mu^\eps_\Lambda \ge
(h_{top} - 4\delta) |\Lambda|
\end{equation}
for any finite open covering $\U^{\eps/2}_\Lambda \in
\C^{\eps/2}_\Lambda$. To finish the proof of the lemma, we use the
sequence of coverings $\{ {\tilde \V}_s\}_s$ defined in Lemma
\ref{lem sub inf} to obtain
\begin{equation} \label{succ-cov-2}
\int_{\F}\ K(f,\eps/2,\Lambda) \ d\mu^\eps_\Lambda = \lim_{s\to
\infty}\ \int_{\F}\ K(f,\tilde \V_s) \ d\mu^\eps_\Lambda \ge
(h_{top} - 4\delta) |\Lambda|
\end{equation}
where the last inequality is given by (\ref{per-U}). \qed

\begin{lem} \label{misura-locale-spaziale}
For any given $\delta>0$ there exists $\eps_0>0$ such that for any
$\eps< \eps_0$ it holds: there exists a probability measure
$\rho^\eps$ on $\F$, invariant with respect to the
$(\zeta_1,\varphi_1)$ action such that
$$K_{\rho^\eps}(\eps/2) := \int_{\F}\
K_{\rho^\eps}(f,\eps/2) \ d\rho^\eps \ge h_{top} - 4\delta$$
\end{lem}

\noindent {\bf Proof.} In Theorem \ref{thm-1} we proved that
$K_{\rho^\eps}(\eps)$, as defined in (\ref{ergodico-p-2}), is not
dependent on the admissible sequence $\La$ of intervals. Hence we
will restrict to the family of intervals $\La=\set{[0,p]}_p$ for
$p\in \N$. For a fixed $q\in \N$, writing $p=tq+r$ with $0\le r<
q$ and using the sub-additivity property proved in Lemma \ref{lem
sub 3}, we have for all $f\in \F$
\begin{equation} \label{spezzare-p}
\begin{array}{ll}
K(f,\eps/2,[0,p]) \le & \sum_{j=0}^{t-2}\
K(\zeta_l(f),\eps/2,[jq,(j+1)q]) + \\[0.3cm]
& + K(f,\eps/2,[0,l]) + K(f,\eps/2,[(t-1)q+l,p]) + t\lo q
\end{array}
\end{equation}
for all $l=0,\dots,q-1$, where $K(f,\eps/2,[0,l])$ and
$K(f,\eps/2,[(t-1)q+l,p])$ are $O(q)$ as shown in
(\ref{limitat-in-L}). Fixed an $\eps<\eps_0$, for all
$p>\lambda_0(\eps)$, where $\lambda_0(\eps)$ is given as in
(\ref{h-eps-lambda}) and in Lemma \ref{misura-locale}, we denote
by $\mu^\eps_p$ the $\varphi_1$-invariant probability measure
associated to $[0,p]$. For all $p>\lambda_0$ we write using
(\ref{spezzare-p})
\begin{equation} \label{spezzare-p-int}
\int_{\F}\ \frac{K(f,\frac{\eps}{2},[0,p])}{p}\ d\mu^\eps_p \le
\int_{\F}\ \sum_{j=0}^{t-2}\
\frac{K(\zeta_l(f),\frac{\eps}{2},[jq,(j+1)q])}{p}\ d\mu^\eps_p +
\frac{t\lo q +o(p)}{p}
\end{equation}
for all $l=0,\dots,q-1$. Moreover for all $j=0,\dots,t-2$ and all
$l=0,\dots,q-1$ it holds
\begin{equation} \label{int-j}
\int_{\F}\ K(\zeta_l(f),\eps/2,[jq,(j+1)q])\ d\mu^\eps_p =
\int_{\F}\ K(f,\eps/2,[0,q])\ d(\zeta^*_{l+jq} \mu^\eps_p)
\end{equation}
hence we define
\begin{equation} \label{media-misure}
\rho^\eps_p:=\frac 1 q \ \sum_{l=0}^{q-1}\ \left( \frac{1}{t-1} \
\sum_{j=0}^{t-2}\ \zeta^*_{l+jq} \mu^\eps_p \right) =
\frac{1}{(t-1)q} \ \sum_{i=0}^{(t-1)q-1}\ \zeta^*_i \mu^\eps_p
\end{equation}
and obtain, using Lemma \ref{misura-locale} and
(\ref{spezzare-p-int})
\begin{equation} \label{fine-spazio}
h_{top}-4\delta \le \frac{(t-1)q}{tq+r} \int_{\F}\
\frac{K(f,\frac{\eps}{2},[0,q])}{q} \ d \rho^\eps_p +\frac{t\lo q
+o(p)}{tq+r}
\end{equation}
We now apply Lemma \ref{convergenza-misure} to the sequence of
measures $\set{\rho^\eps_p}_p$, and we obtain a probability
measure $\rho^\eps$, invariant with respect to the
$(\zeta_1,\varphi_1)$ action, which satisfies
\begin{equation} \label{fine-spazio-2}
\int_{\F}\ \frac{K(f,\frac{\eps}{2},[0,q])}{q} \ d \rho^\eps
+\frac{\lo q}{q} \ge h_{top}-4\delta
\end{equation}
for all $q\in \N$. Letting $q\to \infty$ we have
\begin{equation} \label{fine-spazio-3}
K_{\rho^\eps}(\eps/2) \ge h_{top}-4\delta
\end{equation}
by using the definition of $K_{\rho^\eps}(\eps/2)$ as given in
(\ref{ergodico-p}) and (\ref{ergodico-p-2}). \qed

\vskip 0.5cm In Theorem \ref{thm-1} we have proved that
$K_\mu(\eps)$ is non-decreasing in $\eps$ for any probability
invariant measure $\mu$. Hence, from Lemma
\ref{misura-locale-spaziale}, we obtain that for all
$\eps<\eps_0(\delta)$ it holds
\begin{equation} \label{mis-eps}
K_{\rho^\eps} \ge K_{\rho^\eps}(\eps')\ge h_{top}-4\delta \qquad
\forall\ \eps'< \frac \eps 2
\end{equation}
where we recall that $\rho^\eps$ are probability measure invariant
with respect to the $(\zeta_1,\varphi_1)$ action. To finish the
proof of the theorem, we only need to construct a probability
measure satisfying (\ref{mis-eps}), which is invariant with
respect to space translation and time evolution for all $(x,t)\in
\R \times \R$.

Let us choose a fixed $\eps<\eps_0$ and denote $\rho:=\rho^\eps$.
The probability measure
\begin{equation} \label{misura-totale}
\nu:= \int_0^1  \int_0^1\ \zeta^*_x (\varphi^*_{_{-t}} \rho)\ dt
dx
\end{equation}
is invariant with respect to space translation and time evolution
for all $(x,t)\in \R \times \R$ by definition. We now prove
\begin{lem} \label{fine-princ-var}
The probability measure $\nu$ defined in (\ref{misura-totale})
satisfies
\begin{equation} \label{rho}
K_\nu \ge h_{top}-4\delta
\end{equation}
\end{lem}

\noindent {\bf Proof.} It is enough to prove that for $\eps'$
small enough it holds
\begin{equation} \label{rho-eps}
K_\nu(\eps') \ge h_{top}-4\delta
\end{equation}
For all $p\in \N$ let us write
\begin{equation} \label{fubini-tonelli}
\int_{\F}\ \frac{K(f,\eps',[0,p])}{p} \ d\nu = \int_0^1 \int_0^1
\left( \int_{\F}\ \frac{K(f,\eps',[0,p])}{p}\
d(\zeta^*_x(\varphi^*_{_{-t}} \rho))\right)\ dt dx
\end{equation}
and for the moment consider $(x,t)$ fixed. The first step is to
write
\begin{equation} \label{primo-passo}
\begin{array}{c}
\int_{\F}\ K(f,\eps',[0,p])\ d(\zeta^*_x(\varphi^*_{_{-t}} \rho))
= \int_{\F}\ K(\zeta_x(f),\eps',[0,p])\ d(\varphi^*_{_{-t}}
\rho)=\\[0.5cm]
=\int_{\F}\ K(f,\eps',[x,p+x])\ d(\varphi^*_{_{-t}} \rho)
\end{array}
\end{equation}
By using the sub-additivity property proved in Lemma \ref{lem sub
3}, we write
$$\int_{\F}\ K(f,\eps',[0,p+1])\ d(\varphi^*_{_{-t}} \rho) \le
\int_{\F}\ K(f,\eps',[x,p+x]) \ d(\varphi^*_{_{-t}} \rho) + o(p)$$
where the term $o(p)$ contains the constant $2\lo q$, and the
integral of the terms $K(f,\eps',[0,x])$ and
$K(f,\eps',[p+x,p+1])$ which are bounded as in
(\ref{limitat-in-L}). Hence
\begin{equation} \label{parte-spaziale}
\lim\limits_{p\to \infty} \int_{\F}\ \frac{K(f,\eps',[0,p+1])}{p}\
d(\varphi^*_{_{-t}} \rho) \le \liminf\limits_{p\to \infty}
\int_{\F}\ \frac{K(f,\eps',[x,p+x])}{p} \ d(\varphi^*_{_{-t}}
\rho)
\end{equation}
where on the left hand side we know that the limit exists because
$\varphi^*_{_{-t}} \rho$ is $(\zeta_1,\varphi_1)$-invariant. We
now want to estimate the left hand side. Let us start by writing
$$\int_{\F}\ K(f,\eps',[0,p])\ d(\varphi^*_{_{-t}} \rho) =
\int_{\F}\ K(\varphi_{_{-t}}(f),\eps',[0,p])\ d\rho=$$
$$=\lim_{s\to \infty}\ \lim_{n\to \infty} \int_{\F}\
\frac{K(\varphi_{_{-t}}(f),\tilde \V_s,n)}{n}\ d\rho$$ where we
used the sequence of open coverings $\{ \tilde \V_s \}$ in
$\C^{\eps'}_{[0,p]}$ defined in Lemma \ref{lem sub inf}. By
definition of $K(\varphi_{_{-t}}(f),\tilde \V_s,n)$ we look at the
complexity of the symbolic words in
$\psi(\varphi_{_{-t}}(f),n,\tilde \V_s)$ (see (\ref{compl-n})).
Hence we have
\begin{equation} \label{stima-compl}
K(\varphi_{_{-t}}(f),\tilde \V_s,n) = K(f,\varphi_t(\tilde
\V_s),n)
\end{equation}
where $\varphi_t(\tilde \V_s)$ is a covering of $\F|_{[0,p]}$.
Indeed, for $g\in \F$ there exists $V_j \in \tilde \V_s$ such that
$\varphi_{_{-t}}(g)\in V_j$, because $\varphi$ is invertible.
Hence $g\in \varphi_t(V_j) \in \varphi_t(\tilde \V_s)$. Moreover,
by assumption \textbf{(A3)}, we have that there are constants
$\gamma>0$, $\Gamma >1$ and $C>0$ such that for $p> 2 C
(\eps')^{-1}$ if $d|_{[0,p]}(f_1,f_2)< \eps'$ then
$$d|_{[2C(\eps')^{-1},p-2C(\eps')^{-1}]}(f_1,f_2)< \Gamma e^{\gamma
t}\eps'$$ This implies that for any $t\in [0,1]$, the covering
$\varphi_t(\tilde \V_s)$ is a covering of
$\F|_{[2C(\eps')^{-1},p-2C(\eps')^{-1}]}$ and each of its set is
contained in a ball of radius $\eta = \Gamma e^{\gamma t}\eps'$.
Hence there exists a covering $\tilde \U_s \in
\C^\eta_{[2C(\eps')^{-1},p-2C(\eps')^{-1}]}$ such that
\begin{equation} \label{stima-compl-2}
K(f,\varphi_t(\tilde \V_s),n) \ge K(f, \tilde \U_s,n)
\end{equation}
for all $f\in \F$ and all $n\in \N$. By using (\ref{stima-compl})
and (\ref{stima-compl-2}) we get for all $s\in \N$
$$\lim_{n\to \infty} \int_{\F}\
\frac{K(\varphi_t(f),\tilde \V_s,n)}{n}\ d\rho \ge \lim_{n\to
\infty} \int_{\F}\ \frac{K(f,\tilde \U_s,n)}{n}\ d\rho \ge$$
$$\ge \int_{\F} K(f,\eta,[2C(\eps')^{-1},p-2C(\eps')^{-1}])\ d\rho$$ hence
\begin{equation} \label{parte-tempo}
\begin{array}{c}
\lim\limits_{p\to \infty} \int_{\F}\
\frac{K(f,\eps',[0,p+1])}{p}\ d(\varphi^*_t \rho) \ge \\[0.5cm]
\ge \lim\limits_{p\to \infty} \int_{\F}\
\frac{p+1-4C(\eps')^{-1}}{p}\
\frac{K(f,\eta,[2C(\eps')^{-1},p+1-2C(\eps')^{-1}])}{p+1-4C(\eps')^{-1}}\
d\rho = K_\rho(\eta)
\end{array}
\end{equation}
At this point, we would put together (\ref{primo-passo}),
(\ref{parte-spaziale}) and (\ref{parte-tempo}), and use Lemma
\ref{limsup-K-L} to get
$$K_\nu(\eps') \ge K_\rho(\eta)$$
The assertion would then follow by choosing $\eps'$ small enough
to have $\eta< \frac \eps 2$ and use (\ref{mis-eps}). The only
problem to this argument is that we have tacitly assumed that it
is possible to exchange the order of limit in $p$ and integrations
in $(x,t)$ in (\ref{fubini-tonelli}). However, since by Lemma
\ref{lem sub 3} we have
$$\frac{K(f,\eps',[0,p])}{p} \le K(f,\eps',[0,1]) + \lo q \le const$$
for all $f\in \F$, we can integrate with respect to
$\zeta^*_x(\varphi_{_{-t}} \rho)$ and apply Lebesgue dominated
convergence theorem. \qed

\vskip 0.5cm Hence Theorem \ref{thm-2} is proved.

\appendix

\section{Proof of Lemma \ref{lemma-succ}}

Let us denote as usual
\begin{equation} \label{esse}
(S_n \vartheta)(x):= \sum_{j=0}^{n-1} \vartheta(T^j(x))
\end{equation}
From Birkhoff ergodic theorem we have that there exists a set $Y_1
\subset X$ with $\nu(Y_1^c)=0$ such that for any diverging
sequence $(n_k)_k\subset \Z$  it holds
\begin{equation} \label{birkhoff}
\bar \vartheta(x) := \lim\limits_{k\to \infty}\
\frac{(S_{n_k}\vartheta)(x)}{|n_k|} = \lim\limits_{k\to \infty}\
\frac{1}{|n_k|} \sum\limits_{j=0}^{n_k-1} \vartheta(T^j(x))
\end{equation}
exists, is finite for all $x\in Y_1$ and it is in $L^1(X,\nu)$.
Moreover $\bar \vartheta$ satisfies (\ref{lemma-succ-3}). Let
$Y_1$ be such that the same holds for $|\vartheta|$.

To prove (\ref{lemma-succ-1}), given any sequence of integers
$\set{a_k}_k$ and $\set{b_k}_k$ as in the hypothesis, we write for
all $x\in Y_1$
\begin{equation} \label{vartheta}
\frac{1}{b_k-a_k} \sum\limits_{j=a_k}^{b_k-1} \vartheta(T^j(x)) =
\frac{b_k}{b_k-a_k}\ \frac{(S_{b_k}\vartheta)(x)}{b_k} -
\frac{a_k}{b_k-a_k}\ \frac{(S_{a_k}\vartheta)(x)}{a_k}
\end{equation}
and to study the convergence of (\ref{vartheta}) we divide the
indices $k\in \N$ into four sets:
$$I_1 := \set{k\in \N \ :\ |a_k|\ge \sqrt{b_k-a_k}\ ;\ |b_k|\ge
\sqrt{b_k-a_k}}$$
$$I_2 := \set{k\in \N \ :\ |a_k|\ge \sqrt{b_k-a_k}\ ;\ |b_k|<
\sqrt{b_k-a_k}}$$
$$I_3 := \set{k\in \N \ :\ |a_k|< \sqrt{b_k-a_k}\ ;\ |b_k|\ge
\sqrt{b_k-a_k}}$$
$$I_4 := \set{k\in \N \ :\ |a_k|< \sqrt{b_k-a_k}\ ;\ |b_k|<
\sqrt{b_k-a_k}}$$ First of all we can neglect $I_4$ since it
contains only a finite number of indices by (\ref{cond-succ-1}).
Moreover we introduce for $i=1,2,3$ the notation
$$I_i^{(++)}:= \set{k\in I_i\ :\ a_k\ge 0 \ ;\ b_k\ge 0}$$
and analogously for the other two possible combinations, $(-+)$
and $(--)$. We remark that $I_2= I_2^{(-+)} \cup I_2^{(--)}$ and
$I_3= I_3^{(++)} \cup I_3^{(-+)}$.

Let us consider (\ref{vartheta}) for the indices $k\in
I_1^{(++)}$. By (\ref{cond-succ-1}) and (\ref{birkhoff}), for any
given $x\in Y_1$ and any fixed $\eta >0$ there exists $k_0 \in \N$
such that for all $k\ge k_0$ we have
$$\left| \frac{(S_{b_k}\vartheta)(x)}{b_k} - \bar \vartheta(x)
\right| < \eta$$
$$\left| \frac{(S_{a_k}\vartheta)(x)}{a_k} - \bar \vartheta(x)
\right| < \eta$$ Also by (\ref{cond-succ-2}) we can assume that
for $k\ge k_0$ we have
\begin{eqnarray}
& \frac{a_k}{b_k-a_k} < \limsup\limits_{k\to \infty}
\frac{a_k}{b_k-a_k} + \eta = \frac{1}{l_a}+\eta \label{limsup-a-k}
\\[0.3cm]
& \frac{b_k}{b_k-a_k} < 1+ \limsup\limits_{k\to \infty}
\frac{a_k}{b_k-a_k} + \eta = 1+ \frac{1}{l_a}+\eta
\label{limsup-b-k}
\end{eqnarray}
Applying these inequalities to (\ref{vartheta}) we have that for
all $k\ge k_0$
$$\left| \frac{1}{b_k-a_k} \sum\limits_{j=a_k}^{b_k-1} \vartheta(T^j(x))
- \bar \vartheta(x)  \right| < \eta\ \left(1+ \frac{1}{l_a}+\eta
\right) \left( \frac{1}{l_a}+\eta \right)$$ This proves
(\ref{lemma-succ-1}) for all sequences $a_k$ and $b_k$ with $k\in
I_1^{(++)}$. The same argument applies to $k$ in $I_1^{(-+)}$ and
$I_1^{(--)}$ by writing the right hand side of (\ref{vartheta})
respectively as

\begin{eqnarray}
& \frac{b_k}{b_k+|a_k|}\ \frac{(S_{b_k}\vartheta)(x)}{b_k} +
\frac{|a_k|}{b_k+|a_k|}\ \frac{(S_{a_k}\vartheta)(x)}{|a_k|}
\label{vartheta-neg-pos} \\[0.3cm]
& \frac{|a_k|}{|a_k|-|b_k|}\ \frac{(S_{a_k}\vartheta)(x)}{|a_k|} -
\frac{|b_k|}{|a_k|-|b_k|}\ \frac{(S_{b_k}\vartheta)(x)}{|b_k|}
\label{vartheta-neg-neg}
\end{eqnarray}
and using conditions (\ref{cond-succ-1})-(\ref{cond-succ-3}).

Let us consider now $k\in I_2^{(-+)}$. First of all it holds
\begin{equation} \label{utile-b-k}
0\le \lim_{k\to \infty} \ \frac{b_k}{b_k+|a_k|} < \lim_{k\to
\infty} \ \frac{1}{\sqrt{b_k+|a_k|}} =0
\end{equation}
hence
\begin{equation} \label{utile-a-k}
\lim_{k\to \infty} \ \frac{|a_k|}{b_k+|a_k|} =1
\end{equation}
Moreover we can apply (\ref{birkhoff}) to
$\frac{(S_{a_k}\vartheta)(x)}{|a_k|}$ and, since $b_k <
\sqrt{b_k+|a_k|}$, it holds
\begin{equation} \label{furbizia1}
\limsup\limits_{k\to \infty} \left|
\frac{(S_{b_k}\vartheta)(x)}{b_k+|a_k|} \right| \le
\limsup\limits_{k\to \infty} \frac{1}{\sqrt{b_k+|a_k|}} \
\frac{(S_{_{\sqrt{b_k+|a_k|}}} |\vartheta|)(x)}{\sqrt{b_k+|a_k|}}
=0
\end{equation}
Hence, using (\ref{vartheta-neg-pos}) and applying
(\ref{furbizia1}) to the first term and (\ref{utile-a-k}) and
(\ref{birkhoff}) to the second term, we prove (\ref{lemma-succ-1})
for $k\in I_2^{(-+)}$.

The same arguments apply also to $I_2^{(--)}$, and to $I_3$ by
interchanging the role of $a_k$ and $b_k$.

The proof of (\ref{lemma-succ-2}) follows by a similar argument.
First we show that there exists a set $Y_2\subset X$ with
$\nu(Y_2^c)=0$ such that for any diverging sequence of integers
$\set{n_k}_k$ it holds
\begin{equation} \label{l1-positiva}
\lim\limits_{k\to \infty}\ \frac{\xi(T^{n_k}(x))}{|n_k|} =0
\end{equation}
for all $x\in Y_2$. Since $\xi \in L^1$ and $\xi \ge 0$, and since
the transformation $T$ is measure preserving, for all $\eta >0$ it
holds
$$\sum_{k=1}^\infty\ \nu \set{\xi(T^{k}(x)) > k \eta} =
\sum_{k=1}^\infty \ \nu \set{\xi(x) > k \eta} =$$
$$= \frac{1}{\eta} \sum_{k=1}^\infty \ k \eta \ \nu \set{(k+1)
\eta \ge \xi(x) > k \eta} < \frac{1}{\eta} \sum_{k=1}^\infty
\int_{\set{(k+1) \eta \ge \xi(x) > k \eta}} \ \xi(x) \ d\nu <
$$
$$< \frac{1}{\eta} \ \int_X \ \xi(x) d\nu < \infty$$
hence from the Borel-Cantelli lemma it follows that the measure of
the set on which $\frac{\xi(T^{k}(x))}{k}> \eta$ infinitely often
is zero. Let us moreover assume that the function $\xi(x)$
satisfies Birkhoff theorem (condition (\ref{birkhoff})) for all
$x\in Y_2$.

We now use (\ref{l1-positiva}) as we used (\ref{birkhoff}) before.
If $|b_k|\ge \sqrt{b_k-a_k}$ then by (\ref{cond-succ-1}) $b_k \to
\infty$, hence we can apply (\ref{l1-positiva}) to $b_k$. Hence by
using (\ref{limsup-b-k}), we have that for all $x\in Y_2$ and for
any given $\eta>0$ there exists $k_0(x)$ such that for all $k\ge
k_0(x)$ it holds
$$\frac{\xi(T^{b_k}(x))}{b_k-a_k} \le \eta \left( 1+ \frac{1}{l_a}
+ \eta \right)$$ hence (\ref{lemma-succ-2}) holds in $Y_2$ for
these indices.

If instead $|b_k|< \sqrt{b_k-a_k}$, we can apply (\ref{birkhoff})
by writing
$$\limsup\limits_{k\to \infty} \frac{\xi(T^{b_k}(x))}{b_k-a_k}
\le \limsup\limits_{k\to \infty} \frac{1}{\sqrt{b_k-a_k}} \
\frac{(S_{_{\sqrt{b_k-a_k}}} \xi)(x)}{\sqrt{b_k-a_k}} =0$$ using
$\xi \ge 0$. Hence (\ref{lemma-succ-2}) holds in $Y_2$ also in
this case. This finishes the proof of the lemma by choosing
$Y:=Y_1 \cap Y_2$. \qed

\end{document}